\documentclass[11pt]{article}
\usepackage{url}

\usepackage{graphicx,float,multirow,fullpage,amsmath,amssymb}
\usepackage[colorlinks=true, hyperfigures=false,
           urlcolor=blue, 
           pdfhighlight =/O]{hyperref}
\usepackage{float,multirow,fullpage}

\DeclareMathOperator{\minusone}{-}

\def\eps{\varepsilon}

\newcommand{\eref}[1]{(\ref{#1})}

\newtheorem{theorem}{Theorm}
\newtheorem{prop}{Proposition}
\newtheorem{cor}{Corollary}

\begin{document}

\renewcommand{\thefootnote}{\fnsymbol{footnote}}

\begin{center}
{\bf{\LARGE Mixed-mode oscillations in a multiple time scale \\[0.2cm]
phantom bursting system}}\\[0.5cm]
\end{center}

\vspace{-0.4cm}
\begin{center}
{\bf Maciej KRUPA}$^{1}$\footnotemark[1], \\[0.08cm]
{\bf Alexandre VIDAL}$^{2}$\footnotemark[2], \\[0.08cm]
{\bf Mathieu DESROCHES}$^{1}$\footnotemark[3], \\[0.08cm]
{\bf Fr\'ed\'erique CL\'EMENT}$^{1}$\footnotemark[4]
\end{center}

\vspace{-0.1cm}
\begin{flushleft}
\bf{1} INRIA Paris-Rocquencourt Research Centre, Project-Team SISYPHE, \\
Domaine de Voluceau Rocquencourt - B.P. 105, 78153 Le Chesnay cedex, France.
\\[0.2cm]
\bf{2} Universit\'e d'\'Evry-Val-d'Essonne, Laboratoire Analyse et Probabilit\'es EA 2172, F\'ed\'eration de Math\'ematiques FR 3409, IBGBI, 23 Boulevard de France, 91037, Evry, France.
\\[0.2cm]
\end{flushleft}

\renewcommand{\thefootnote}{\fnsymbol{footnote}}
\setcounter{footnote}{1}
\footnotetext[1]{{\tt Maciej.P.Krupa@gmail.com}}
\footnotetext[2]{\tt Alexandre.Vidal@univ-evry.fr}
\footnotetext[3]{\tt Mathieu.Desroches@inria.fr}
\footnotetext[4]{\tt Frederique.Clement@inria.fr}
\renewcommand{\thefootnote}{\arabic{footnote}}

\vspace{-0.4cm}
\section*{Abstract}
In this work we study mixed mode oscillations in a model of secretion of GnRH (Gonadotropin Releasing Hormone). The model is a phantom burster consisting of two feedforward coupled FitzHugh-Nagumo systems, with three time scales. The forcing system (Regulator) evolves on the slowest scale and acts by moving the slow null-cline of the forced  system  (Secretor). There are three modes of dynamics: pulsatility (transient relaxation oscillation), surge (quasi steady state) and small oscillations related to the passage of the slow null-cline through a fold point of the fast null-cline. We derive a variety of reductions, taking advantage of the mentioned features of the system. We obtain two results; one on the local dynamics near the fold in the parameter regime corresponding to the presence of small oscillations and the other on the global dynamics, more specifically on  the existence of an attracting limit cycle. Our local result is a rigorous characterization of small canards and sectors of rotation in the case of folded node with an additional time scale, a feature allowing for a clear geometric argument. The global result gives the existence of an attracting unique limit cycle, which, in some parameter regimes, remains attracting and unique even during passages through a canard explosion.

\vspace{0.4cm}
\noindent \textbf{Keywords :} Slow-fast systems, multiple time scales, mixed mode oscillations, limit cycles, se\-condary canards, sectors of rotations, folded node, singular perturbation, blow-up, GnRH secretion. \\

\vspace{0.01cm}
\noindent \textbf{AMS Classification :} \\[0.2cm]
\begin{tabular}{ll}
34C15 Nonlinear oscillations, coupled oscillators, & 34C23 Bifurcation, \\
34C26 Relaxation oscillations, & 34D15 Singular perturbations, \\
34E13 Multiple scale methods, & 34E15 Singular perturbations, general theory, \\
34E17 Canard solutions, & 70K70 Systems with slow and fast motions, \\
92B05 General biology and biomathematics. & \\
\end{tabular}

\section{Introduction}
Mixed Mode Oscillations (MMOs) is a term used to describe trajectories that combine small oscillations and large oscillations of relaxation type, both recurring in an alternating manner. Recently there has been a lot of interest in MMOs that arise due to a generalized canard phenomenon, starting with the work of Milik, Szmolyan, Loeffelmann and Groeller \cite{am-ps-hl-eg_98}. Such MMOs arise in the context of slow-fast systems with at least two slow variables and with a folded critical manifold (set of equilibria of the fast system). The small oscillations arise during the passage of the trajectories near a fold, due to the presence of a so-called \textit{folded singularity}. The dynamics near the folded singularity is transient, yet recurrent: the trajectories return to the neighborhood of the folded singularity by way of a global return mechanism. 

An important step on the way to an understanding of MMOs is the analysis of the flow near the folded singularities. Of particular importance are special solutions called {\it canards}. The term canard was first used to denote periodic solutions of the van der Pol equation that stayed close to the unstable slow manifold (approximated by the middle branch of the fast nullcline) \cite{eb-jlc-fd-md_81}. One of the characteristic features of canard cycles is that they exist only for an exponentially small range of parameter values. This very sharp transition was then termed {\it canard explosion} \cite{mb_88}. The related term canard solution has been used to denote solutions connecting from a stable slow manifold to an unstable slow manifold. Such canards, sometimes also called {\it maximal canards}, organize the dynamics in a similar way as invariant sets which separate different dynamical regimes (e.g., separatrices of saddle points). In systems with more than one slow variable, canards occur in a more robust fashion and underlie the presence of the small oscillations near the folded singularity in MMOs.

A prototypical example of a folded singularity with small oscillations is the \textit{folded node}, studied by Beno\^it \cite{eb_90}, by Wechselberger and Szmolyan \cite{ps-mw_01}, and by Wechselberger \cite{mw_05}. These articles focused on the local aspects of the dynamics. An exposition of how the dynamics near the folded node can be combined with a global return mechanism to lead to MMOs was given in \cite{mb-mk-mw_06}. This work was used as a basis of various explanations of MMO dynamics found in applications \cite{rw_07, hr-mw-nk_08, be-mw_09}. A shortcoming of the folded node approach is the lack of connection to a Hopf bifurcation, which seems to play a prominent role in many MMOs.  This led to the interest in another, more degenerate folded singularity, known as \textit{Folded Saddle Node of type II} (FSNII), originally introduced in  \cite{am-ps-hl-eg_98} and recently analyzed in some detail by Krupa and Wechselberger \cite{mk-mw_10}. Guckenheimer \cite{jg_08} studied a very similar problem in the parameter regime yet closer to the Hopf bifurcation, calling it \textit{singular Hopf bifurcation}. For a more comprehensive overview we refer the reader to the recent review article \cite{md-jg-bk-ck-ho-mw_12}.

Two notions that are central to the study of MMOs are {\it secondary canards} and {\it sectors of rotation}. Secondary canards \cite{mw_05, mb-mk-mw_06} are trajectories which originate in the attracting slow manifold, make a number of small oscillations in the fold region, and continue to the unstable slow manifold. There is ample numerical evidence of the existence and role of secondary canards (\cite{md-bk-ho_08,md-bk-ho_08_2}), as well as some partial theoretical results (\cite{jg-rh_05,mw_05, mk-mw_10, mk-np-nk_08}). It has been highlighted that two trajectories crossing the region between two consecutive canards display the same number of small oscillations. Hence, the regions separated by secondary canards have been called  {\it sectors of rotation} (\cite{mb-mk-mw_06}). As a parameter changes, a periodic orbit may move closer to a canard  and pass to the adjacent sector of rotation. This transition has never been studied in detail. It is similar to a canard explosion, although more complicated, as chaotic behavior can be expected.

The main result of this paper is that, in the context of our phantom burster problem, there exists an attracting MMO orbit for all parameter values, also during the passage between different sectors. In addition we obtain a result on the existence of secondary canards of rotational type that is complementary to the results in \cite{mw_05, mk-np-nk_08, mk-mw_10} and relevant to the context of the phantom burster.

It is important to note that, even if canards are more robust in three-dimensional slow-fast systems, they are still difficult to find numerically as well as particular types of MMOs. Forward integration is not possible due to the exponential expansion along the slow manifold, leading to exponential magnification of numerical errors. A breakthrough in the numerical detection and continuation of canards and MMOs has been achieved by Desroches et. al. \cite{md-bk-ho_08} who used a boundary value approach in the context of numerical continuation with the software package \textsc{Auto}.
 
In this article we investigate the presence of MMOs in the following system: 
\begin{subequations} \label{GnRHSystem}
\begin{eqnarray} 
\eps \delta \dot{x} &=& -y+f(x), \label{eqdotx0}\\
\delta \dot{y} &=& a_0 x+a_1 y+a_2+c X, \label{eqdoty0}\\
\delta \dot{X} &=& -Y+g(X), \label{eqdotX0}\\
\dot{Y} &=& X+b_1Y+b_2, \label{eqdotY0}
\end{eqnarray}
\end{subequations}
with
\[
\begin{array}{c}
f:x\rightarrow \lambda_3 x^3+ \lambda_1 x, \\
g:x\rightarrow \mu_3 x^3+ \mu_1 x, \\
\lambda_3 , \mu_3 <0, \quad \lambda_1, \mu_1 > 0, \\
a_i, c>0,\\
0<\eps , \delta <<1.
\end{array}
\]
System \eqref{GnRHSystem} has been proposed in \cite{fc-jpf_07, fc-av_09, av-fc_10} to model the dynamics of GnRH secretion by hypothalamic neurons in female mammals. Subsystem \eref{eqdotx0}-\eref{eqdoty0}, called the {\it Secretor}, represents the mean-field approximation of the GnRH neuron population dynamics. It is driven, through the coupling term $cX$, by subsystem \eref{eqdotX0}-\eref{eqdotY0}, called the {\it Regulator}, representing the activity of the interneuron population that conveys the periodic action of the ovarian steroids onto the GnRH neuron population. The time scale difference between the two oscillators is a transcription of the ratio between the ovarian cycle duration (few weeks) and the period of the secretory activity of the GnRH neuron population (few hours).

From the point of view of dynamical classification, system \eref{GnRHSystem} is a {\it phantom burster} with the additional feature that it has multiple time scales. Both the Regulator \eref{eqdotX0}-\eref{eqdotY0} and the driven Secretor \eref{eqdotx0}-\eref{eqdoty0} are slow-fast systems of Fitzhugh-Nagumo type and there is a time scale difference between the two. The parameters controlling the time scales are $\eps $ and $\delta $.

\begin{figure}[htbp]
\centering
\includegraphics[scale=0.68]{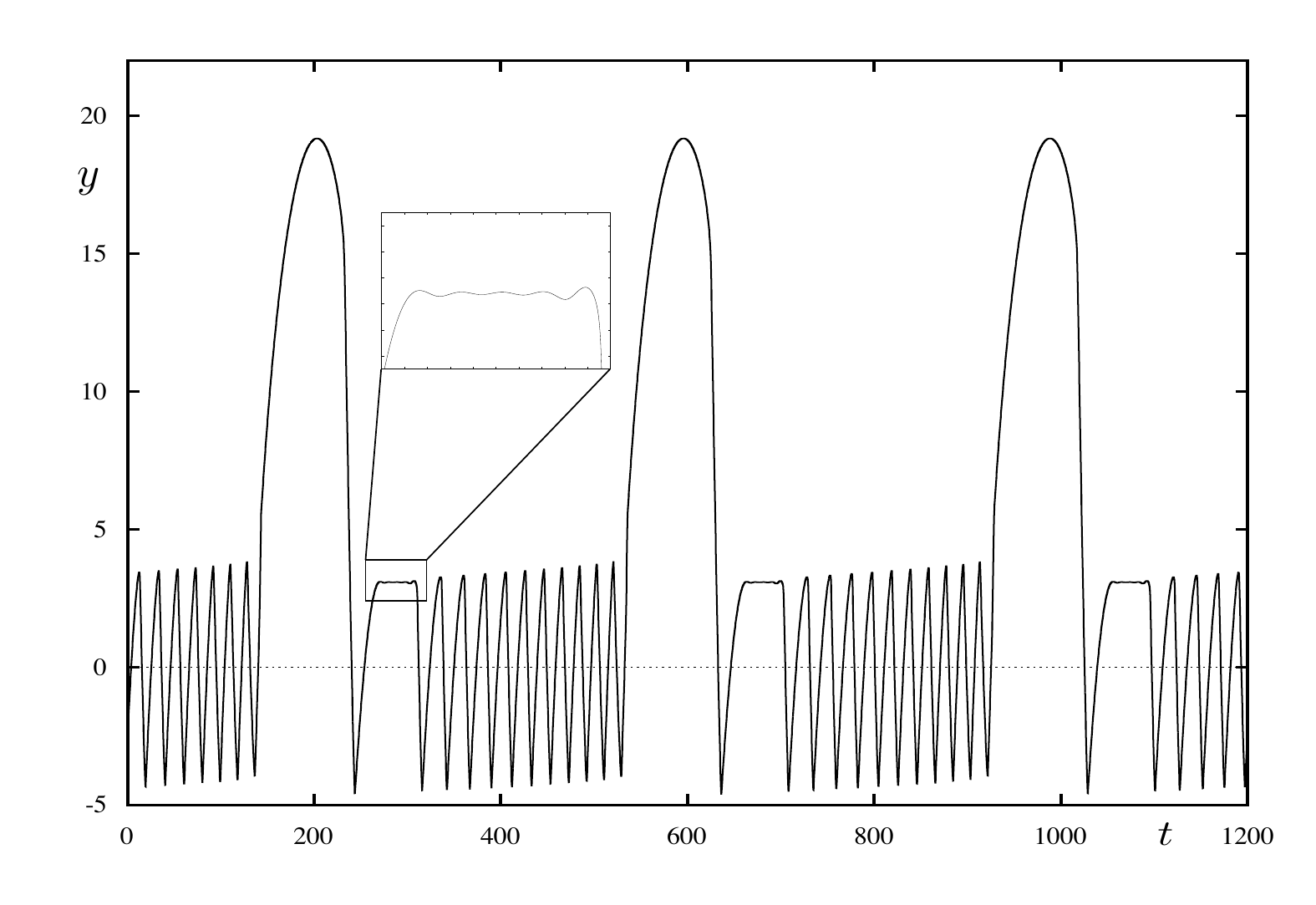}
\caption{The signal $y(t)$ generated along a typical orbit of system \eref{GnRHSystem} displays a periodic alternation of pulsatile phases, surges and pauses. As displayed in the inset, the pauses consist of small oscillations.}
\label{GnRHSecr}
\end{figure}

We are interested in the case when the Regulator displays a stable relaxation limit cycle. Then, in a certain region of the parameter space (see \cite{fc-av_09}), the $X$-driven Secretor alternates between a fast oscillatory regime (when it displays an attracting relaxation limit cycle) and a stationary regime (when the current point tracks an attracting singular point). As a result, the signal generated by $y$, that represents the GnRH secretion along time, displays a periodic alternation of pulsatile regimes and surges as illustrated in Figure \ref{GnRHSecr}. During the transition from a surge back to the subsequent pulse phase, a {\it pause} corresponding to a segment of small oscillations may occur.

In this article, we analyze the dynamical mechanism based on three different time scales that underlies the occurrence of the small oscillations. We prove that, for certain choices of the parameter values, the MMOs, including the pulse phase, surge and pause, exist and are stable limit cycles, even when close to a secondary canard. More precisely, we prove that canards with a specified number of small oscillations are unique (with fixed choices of slow manifolds) and that any two adjacent canards differ by one rotation. Thus we prove that sectors of the same rotation (or simply sectors of rotation) exist and the passage, as a parameter varies, through a secondary canards adds (or subtracts) one small oscillation to the globally attracting orbit.

The paper is organized as follows. In \S \ref{sec-pbi} we present the phantom burster dynamics of \eqref{GnRHSystem}, discuss the different phases of the orbits and state the main result (Theorem \ref{th-main}). Section \ref{sec-mr} is devoted to the local analysis of canard oscillations in a three time scale reduced system (Theorem \ref{th-seca}). In \S \ref{sec-global} we analyze the return mechanism and prove Theorem \ref{th-main}. Section \ref{sec-num} contains numerical findings illustrating our results. The article ends with a discussion section.

\section{Multiple time scale phantom bursting}\label{sec-pbi}
In this section we give a qualitative description of the dynamics we are interested in and state the main result. Parts of this section are a review and we refer the reader to \cite{fc-jpf_07, fc-av_09} for more details.  We begin by sketching the basic features of the dynamics of the decoupled Regulator \eqref{eqdotX0}-\eqref{eqdotY0}. Subsequently  we set constraints on the Secretor's parameters in order to obtain the right dynamical behavior, introduce the different reduced systems suitable to describe various stages of the dynamics, and  briefly describe the evolution of the $x$ and $y$ variables in the different stages of the dynamics. We end the section by stating the main theorem.

\subsection{The Regulator dynamics and its influence on the position of the Secretor slow nullcline and singular points} \label{sec-sketch}

\begin{figure}[htbp]
\centering
\includegraphics[width=0.5\textwidth]{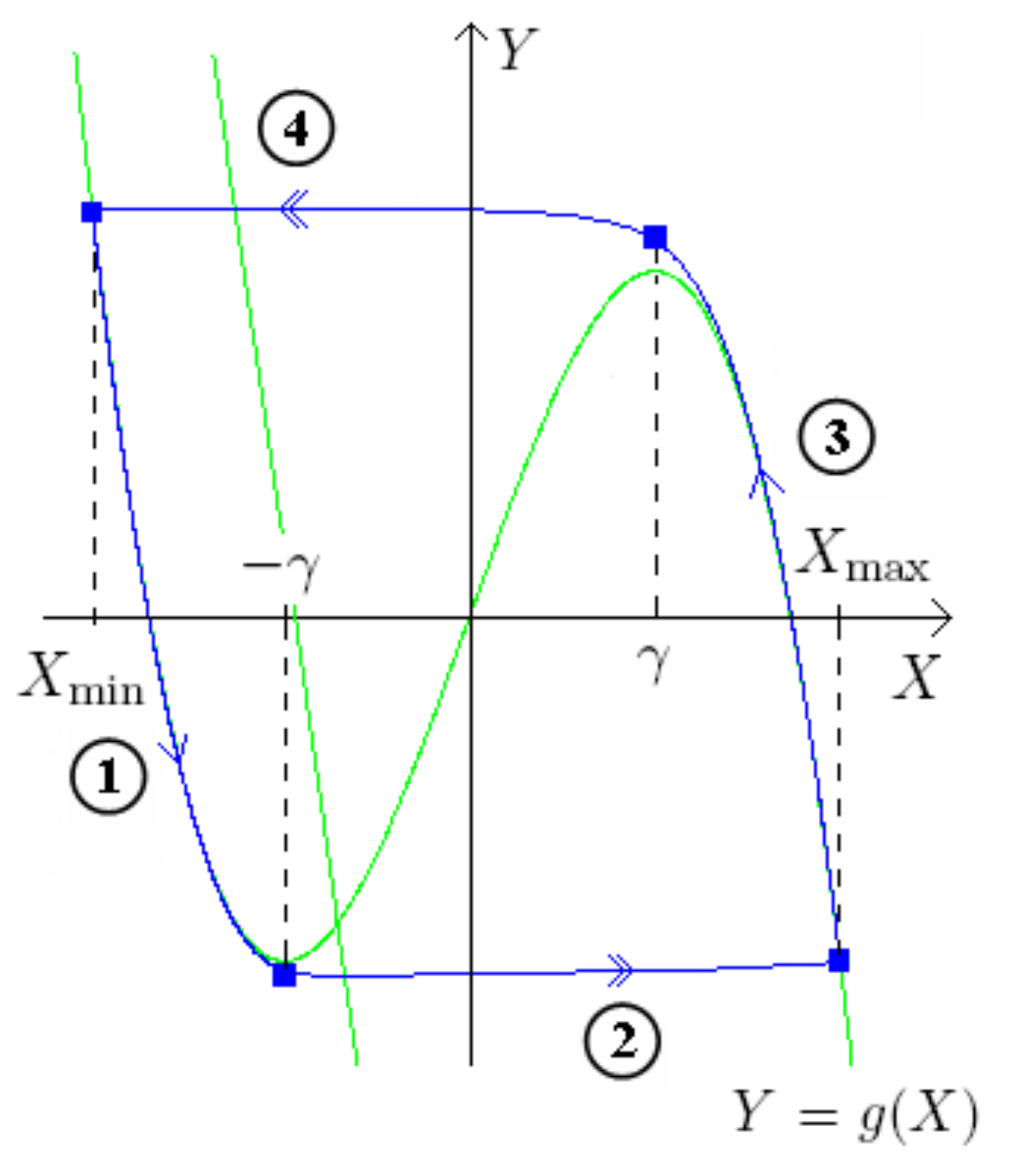}
\caption{Limit cycle of \eref{eqdotX0}-\eref{eqdotY0} and its four different phases.}
\label{CycleXY}
\end{figure}

We define $\gamma >0$ by $g'(\pm \gamma)=0$ so that the two knees of the cubic $X$-nullcline $Y=g(X)$ are $(\pm \gamma , g(\pm \gamma))$. The knees split the cubic $X$-nullcline  into three parts : the left and right branches where $g'<0$ and the middle where $g'>0$. As mentioned in the introduction, we assume that the parameters, specifically $b_1$ and $b_2$, are chosen so that the Regulator admits a relaxation limit cycle. To ensure this property, it is sufficient to assume that $b_1$ is small enough (so that the $Y$-nullcline is steep enough) and that the $Y$-nullcline intersects the cubic $X$-nullcline $Y=g(X)$ on its middle branch (where $g'<0$) away from the knees (see \cite{fc-av_09}). Let us note $X_{\min}$ and $X_{\max}$ respectively the minimal and maximal value of $X$ along the Regulator limit cycle. For later reference we list four different phases of the evolution of $X$ (see Figure \ref{CycleXY}):
\begin{enumerate}
\item slow motion near the left branch of the cubic $Y=g(X)$ : $X$ increases slowly from $X_{\min}$ to $-\gamma $,
\item fast motion from the left knee to the right branch of the cubic : $X$ increases quickly from $-\gamma $ to $X_{\max}$,
\item slow motion near the right branch of the cubic : $X$ decreases slowly from $X_{\max}$ to $\gamma$,
\item fast motion from the right knee to the left branch of the cubic : $X$ decreases quickly from $\gamma $ to $X_{\min}$.
\end{enumerate}

The value of $X$ drives the Secretor $y$-nullcline  defined by $a_0 x+a_1 y+a_2+c X=0$. Note that this nullcline is a straight line whose slope $-a_0/a_1$ does not depend on $X$. As usual in the Fitzhugh-Nagumo system, $a_1$ is assumed to be small so that the $y$-nullcline is very steep.

As $X$ increases (resp. decreases), the $y$-nullcline moves to the left (resp. to the right) in the Secretor phase space $(x,y)$. Hence, depending on the value of $X$, the number of singular points (lying on the cubic $x$-nullcline $y=f(x)$) varies. Also, their nature depends on their position with respect to the fold points $(x_f,f(x_f))$ and $(-x_f,-f(x_f))$ that splits the $x$-nullcline into three part (left, middle and right branch). In particular:
\begin{enumerate}
\item[a.] if a singular point lies on the middle branch outside a $O(\eps)$-neighborhood of the folds, it is surrounded by a relaxation limit cycle ;
\item[b.] if two different singular points lie on the left (resp. right) branch, the lowest (resp. highest) one is an attracting node and the highest (resp. lowest) one is a saddle ;
\item[c.] if the Secretor admits a unique singular point, it is a saddle.
\end{enumerate}
The passage from a to b is a Hopf bifurcation that makes the limit cycle disappear through a canard explosion and the passage from b to c is a saddle-node bifurcation that makes the saddle and the node collapse.

In the following, we assume (see hypotheses H1 to H4 in the following section) that for all values of $X$ between $X_{\min}$ and $X_{\max}$, the Secretor admits three different singular points determined by their $x$-component. Of special importance is the middle singular point (corresponding to the $x$-value lying between the two others) for which we note the $x$-component $x_{\rm sing}(X)$.

\subsection{Constraints on the Secretor parameters and statement of the main result}

\begin{figure}[htbp]
\centering
\begin{tabular}{ccc}
\hspace{0.3cm} {\bf H1 : $X=X_{\min}$} & \hspace{0.3cm} {\bf H2 : $X=-\gamma$} \\
\includegraphics[width=7cm]{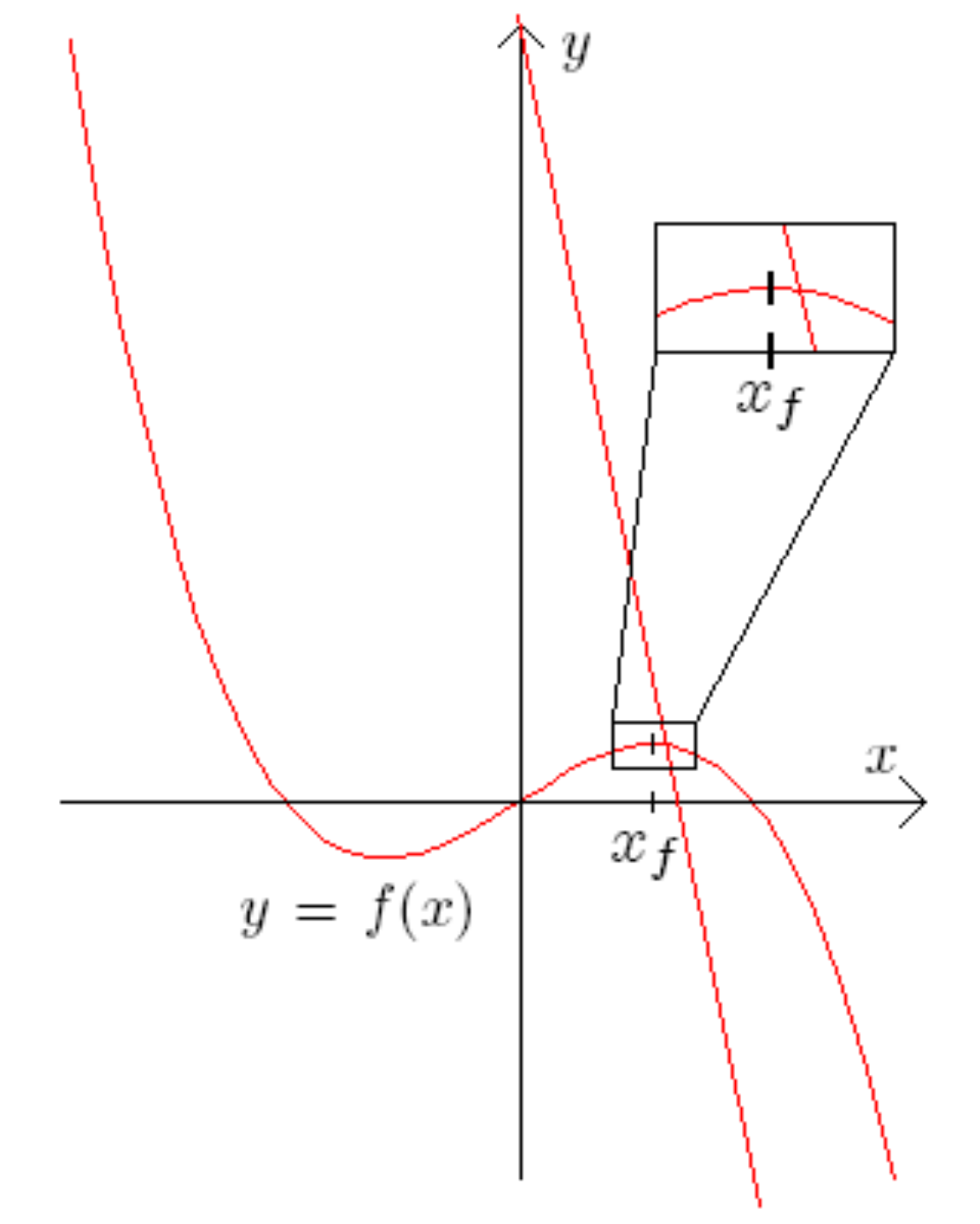} & \includegraphics[width=7cm]{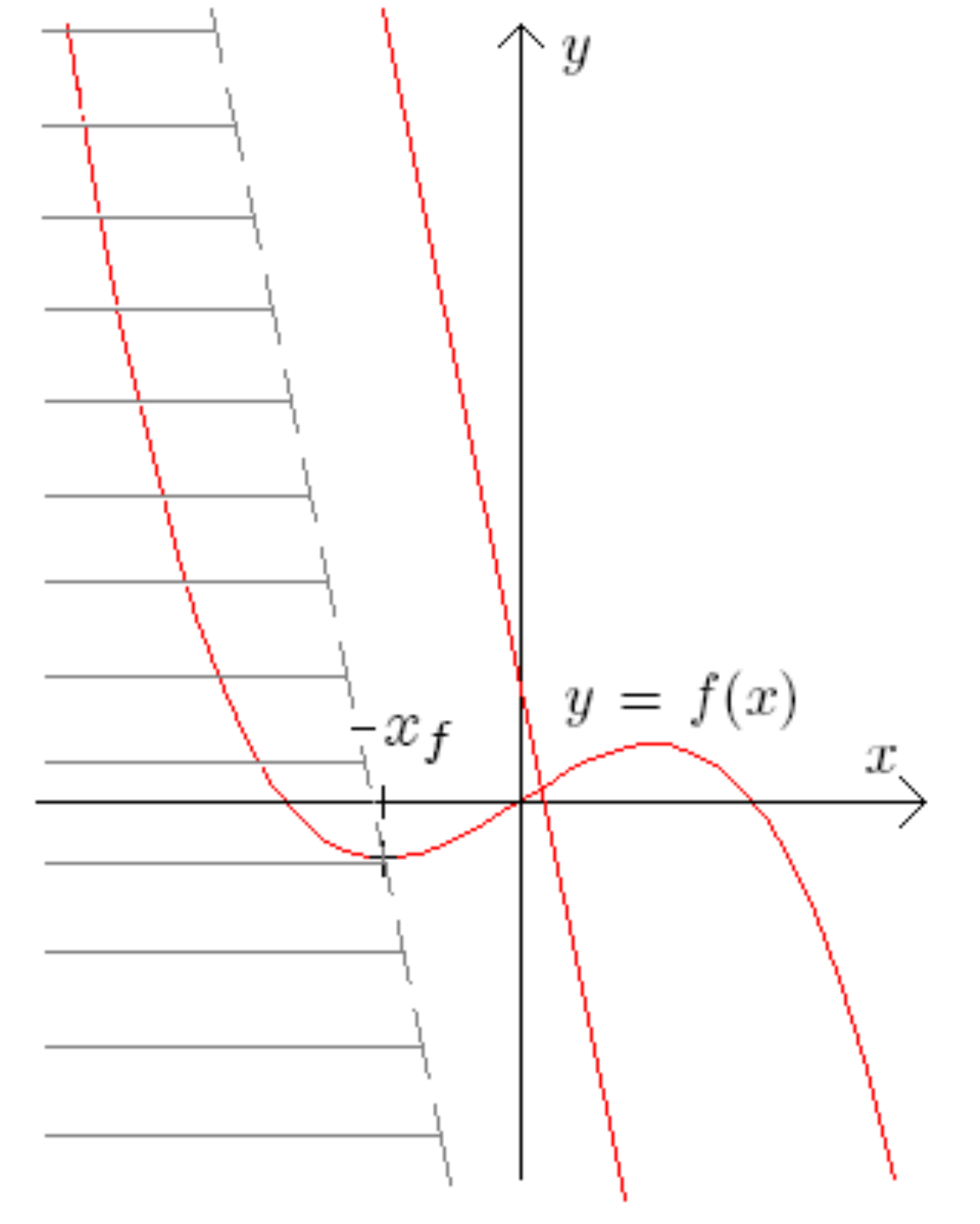} \\
\hspace{0.3cm} {\bf H3 : $X=X_{\max}$} & \hspace{0.3cm} {\bf H4 : $X=\gamma$} \\
\includegraphics[width=7cm]{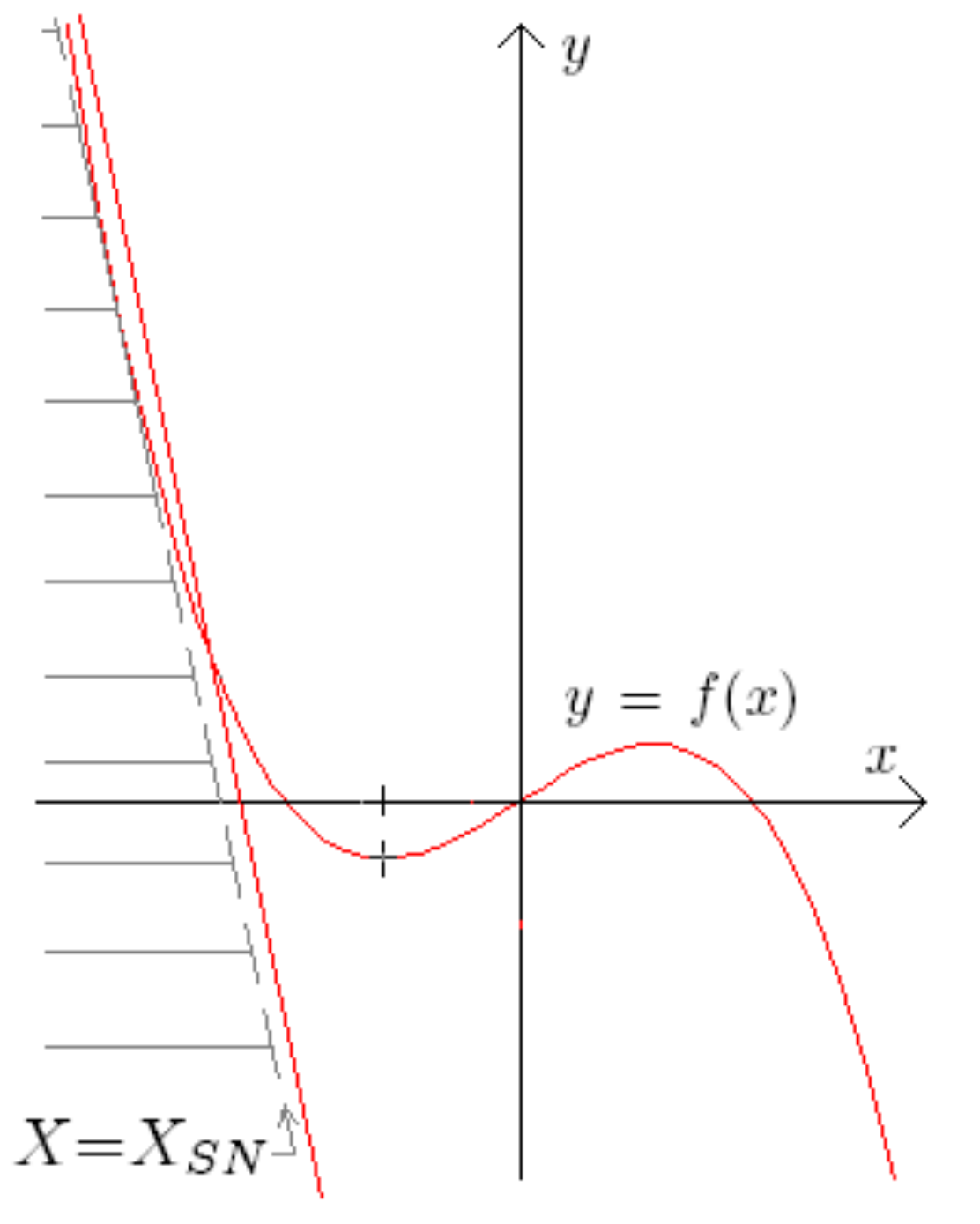} & \includegraphics[width=7cm]{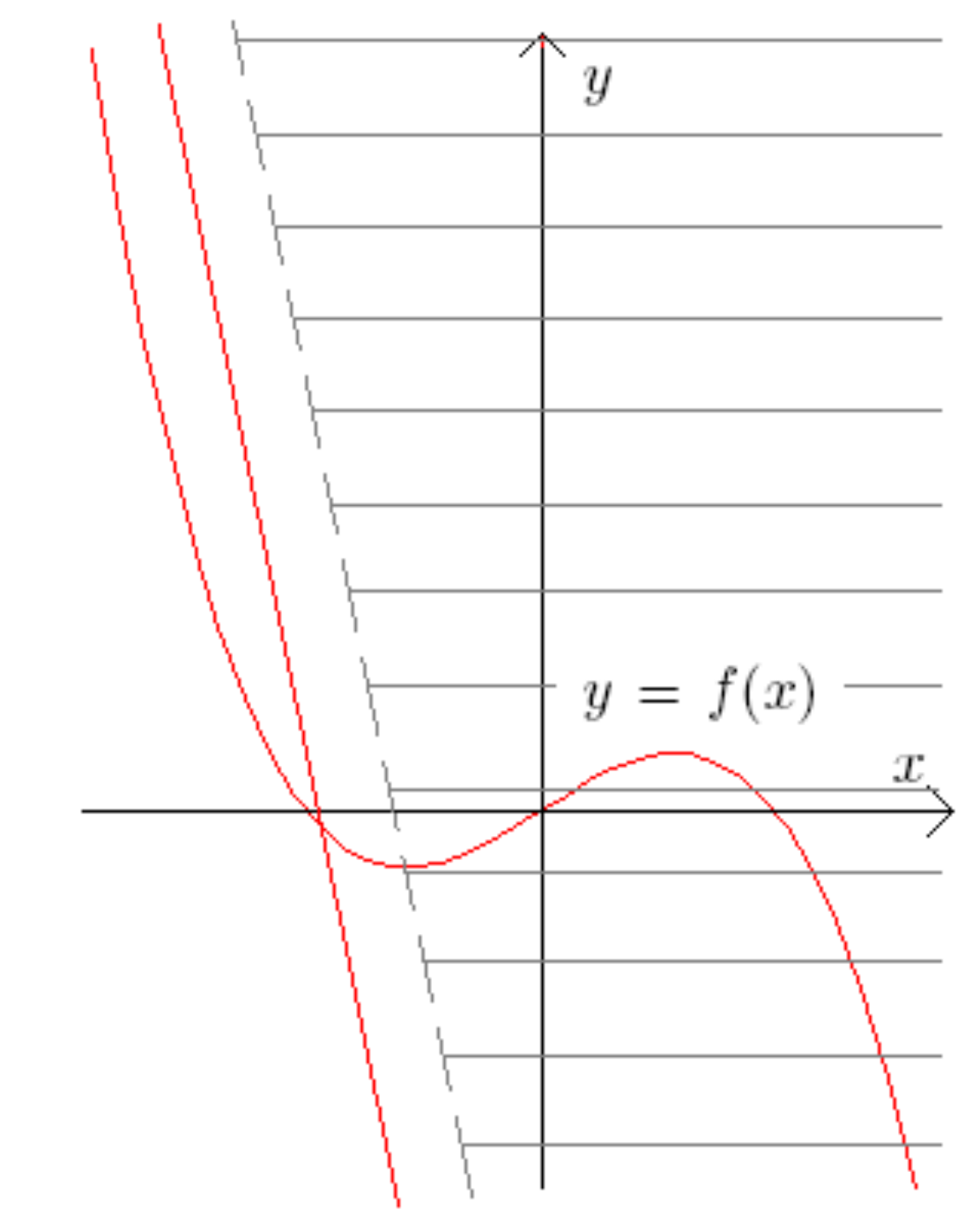}
\end{tabular}
\caption{Illustration of the four hypotheses (H1) to (H4) on parameters to obtain the right system behavior.}
\label{Conds}
\end{figure}

To obtain the qualitative behavior of the $y$-signal generated by full system \eqref{GnRHSystem} (Figure \ref{GnRHSecr}) we make the following hypotheses illustrated by Figure \ref{Conds}.

\begin{description}
\item[(H1)] The $y$-nullcline should pass through the right fold point of the cubic $y=f(x)$ which generates the small oscillations. Hence, we assume that, for $X=X_{\min}$, the $y$-nullcline should be on the right of -- and close to -- the upper fold $(x_f,f(x_f)$:
\[
x_f \lesssim x_{\rm sing}(X_{\min}) \text{ i.e. } X_{\min} \lesssim X_f =-\frac{a_0 x_f + a_1 f(x_f) +a_2}{c}.
\]
This condition will be discussed in more detail later on.
\item[(H2)] Once the $y$-nullcline has passed the right fold and the relaxation limit cycle of system \eref{eqdotx0}-\eref{eqdoty0} appears, the cycle should persist until $X=-\gamma $. Hence, we assume that for $X=-\gamma $, the $y$-nullcline intersect the cubic $y=f(x)$ on its middle branch:
\[
-x_f < x_{\rm sing}(-\gamma) \text{ i.e. } -a_0 x_f+a_1 f(-x_f)+a_2-c\gamma <0.
\]
\item[(H3)] From the beginning of the surge phase, system \eref{eqdotx0}-\eref{eqdoty0} must admit an attracting node and a saddle on the left branch of the cubic $y=f(x)$. This condition reads
\[
a_0+a_1 \lambda_1 > 2 \sqrt{a_1 c \lambda_1 X_{\max}}
\]
which is equivalent to:
\[
X_{\max}<X_{SN}=\frac{(a_0+a_1 \lambda_1)^2}{4a_1 c \lambda_1}.
\]
Let us note that value $X_{SN}$ of $X$ corresponds to the saddle-node bifurcation of the Secretor occurring when the $y$-nullcline is tangent to the left branch of the cubic $y=f(x)$.
\item[(H4)] Until the end of the surge phase, system \eref{eqdotx0}-\eref{eqdoty0} must admit an attracting node and a saddle on the left branch of the cubic $y=f(x)$ as well. This condition reads:
\[
-x_f>x_{\rm sing}(\gamma)  \text{ i.e. } -a_0 x_f+a_1 f(-x_f)+a_2+c\gamma >0.
\]
\end{description}

\begin{figure}[H]
\centering
\includegraphics[width=0.6\textwidth]{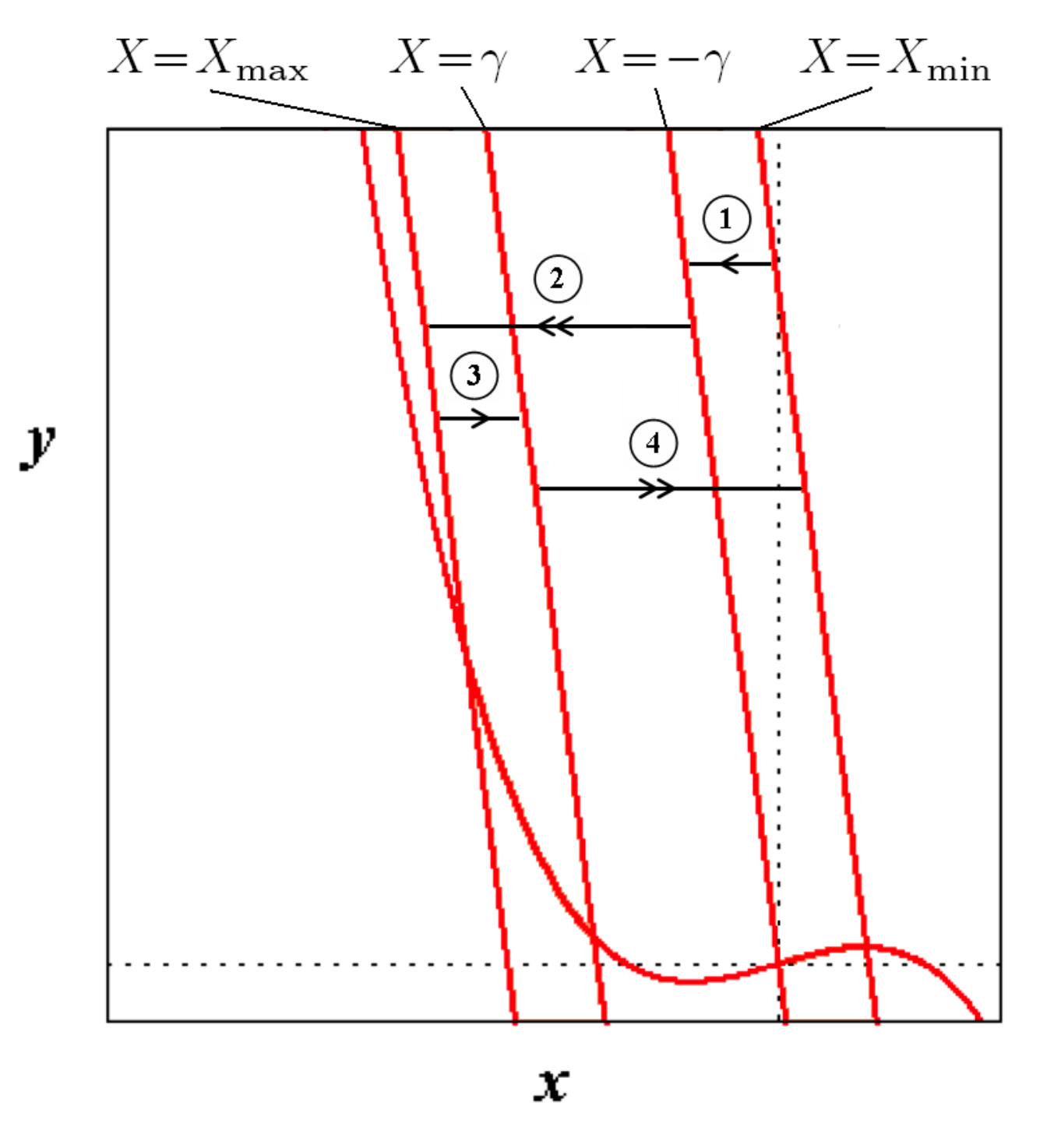}
\caption{Instance of the Secretor slow nullcline locations under hypotheses (H1) to (H4).}
\label{Type1}
\end{figure}

Figure \ref{Type1} gives an instance of the relative positions of the Secretor slow nullcline that can arise due to the variation in $X$ as $(X,Y)$ traces the Regulator relaxation cycle under assumptions H1 to H4. Then, the signal generated by variable $y$ displays an alternation of surge, small oscillations and surge phases as illustrated in Figure \ref{GnRHSecr}. We refer to \cite{av-fc_10} for an expanded explanation of the model behavior based on this approach.

\subsection{Reduced systems} \label{sec-red-sys}

For each phase from 1 to 4 of the Regulator limit cycle, system \eref{GnRHSystem} can be reduced using a specific approximation. We first recall the general process of desingularization applied to slow-fast systems near a fold. We introduce the different reduced systems that we will use for the dynamics analysis.\\

\noindent \textbf{Desingularized Reduced System.} When a general slow-fast dynamical system
\begin{subequations} \label{gen}
\begin{eqnarray} 
\eps\dot{x} &=& f(x,y,\eps), \label{fastgen} \\
\dot{y} &=& g(x,y,\eps), \label{slowgen}
\end{eqnarray}
\end{subequations}
is considered, with $x$ and $y$ of arbitrary dimension and $0<\eps\ll 1$, one classic way to understand the overall dynamics is by looking at the slow and the fast dynamics separately. An object of great importance for both the slow and the fast dynamics of the full system is the so-called {\it critical manifold} $C^0$ defined as the nullcline for the fast variable, that is:
\[
C^0 = \Big\{(x,y);\;f(x,y,0)=0\Big\}.
\]
Consequently, the critical manifold is the phase space of the reduced system obtained by setting $\eps=0$ in equations \eref{gen} and which approximates the slow dynamics of the original system; the reduced system is a differential-algebraic equation. In order to understand the flow of the reduced system, which then takes place on $C^0$ and is associated with the singular limit $\eps=0$, the usual strategy -- which we will use several times in the rest of the paper -- is to differentiate the algebraic equation defining $C^0$ with respect to time. This gives
\begin{eqnarray*}
\dot{x}f_x(x,y,0)+\dot{y}f_y(x,y,0) &=& 0,\\
\dot{y} &=& g(x,y,0),
\end{eqnarray*}
which reduces to 
\begin{subequations} \label{red}
\begin{eqnarray}
\dot{x}f_x(x,y,0) &=& -g(x,y,0)f_y(x,y,0), \label{redx}\\
\dot{y} &=& g(x,y,0). \label{redy}
\end{eqnarray}
\end{subequations}
The previous system is singular along the fold set of $C^0$ with respect to the fast variable $x$, that is, the set $\mathcal{F}=\{f_x(x,y,0)=0\}$. In order to understand the slow flow up to the fold set, one can desingularize system \eref{red} via a rescaling of factor $f_x(x,y,0)$, which yields the so-called {\it desingularized reduced system}. In that case, special care has to be taken going from the desingularized reduced system back to the slow system. Indeed the previous rescaling changes the orientation of orbits when $f_x(x,y,0)<0$ and one needs to reverse orientation in order to get the correct direction of the flow in the original reduced system. \\

\noindent \textbf{Three dimensional reduction with three time scales during the pulsatile phase.}
Under the preceding assumptions, slow motion 1 ($X_{\min} < X < \gamma$) corresponds for system \eref{eqdotx0}-\eref{eqdoty0} to the oscillatory phase producing the small oscillations and subsequently the pulses in the $y$-signal.

The variables $X$ and $Y$ follow the slowest time scale and the current point $(X,Y)$ remains in a $O(\delta)$-neighborhood of the cubic. This reads $Y=h_{\delta}(X)$ where $(X,\delta) \mapsto h_{\delta }(X)$ is an analytic function on $]-\infty, -\gamma[ \times \mathbb{R_+^*}$ and $h_{0}=g$. Thus, on $]-\infty, -\gamma[$, $h^{\prime }_{\delta}(X)=g'(X)+O(\delta)$.

We introduce a reduced system obtained from \eref{GnRHSystem} assuming that $Y=h_{\delta } (X)$. We differentiate this condition, with $\delta $ constant:
\[
\dot{Y}=\dot{X} h^{\prime }_{\delta}(X).
\]
By replacing the dynamics of $\dot{Y}$ in \eref{GnRHSystem}, one obtains  the three-dimensional system with three different time scales:
\begin{subequations} \label{3D3TS}
\begin{eqnarray} 
\eps \delta \dot{x} &=& -y+f(x), \label{3D3TSx}\\
\delta \dot{y} &=& a_0 x+a_1 y+a_2+c X, \label{3D3TSy}\\
\dot{X} &=& \frac{X+b_1 (g(X)+O(\delta)) +b_2}{g'(X)+O(\delta)}. \label{3D3TSX} 
\end{eqnarray}
\end{subequations}

\noindent \textbf{Two-dimensional reduction with two time scales during the surge phase.}
Slow motion 3 ($\gamma < X < X_{\max}$) corresponds to the surge phase. The current point $(x,y)$ follows the attracting node of \eref{eqdotx0}-\eref{eqdoty0} lying on the left branch of $y=f(x)$. Hence, both approximation $X\simeq g(X)$ and $y\simeq f(x)$ stand.

By reducing the fastest time scale, i.e., by setting $y=f(x)$ in \eqref{GnRHSystem}, we obtain the following system:
\begin{subequations} \label{GnRHlayer}
\begin{eqnarray} 
\delta f'(x)\dot{x} &=& a_0 x+a_1 f(x)+a_2+c X, \label{GnRHlayerx}\\
\delta \dot{X} &=& -Y+g(X), \label{GnRHlayerX}\\
\dot{Y} &=& X+b_1Y+b_2. \label{GnRHlayerY}
\end{eqnarray}
\end{subequations}
Then setting $Y=g(X)$ leads to the following equations:
\begin{subequations} \label{GnRHlayer1}
\begin{eqnarray} 
\delta f'(x)\dot{x} &=& a_0 x+a_1 f(x)+a_2+c X, \label{GnRHlayer1x}\\
g'(X)\dot{X} &=& X+b_1g(X)+b_2. \label{GnRHlayer1X}
\end{eqnarray}
\end{subequations}
Away from the folds of both cubics (where $f'(x)=0$ or $g'(X)=0$) we can rewrite \eqref{GnRHlayer1} as follows:
\begin{subequations} \label{GnRHlayer2}
\begin{eqnarray} 
\delta\dot{x} &=& \frac{a_0 x+a_1 f(x)+a_2+c X}{f' (x)}, \label{GnRHlayer2x}\\
\dot{X} &=& \frac{X+b_1g(X)+b_2}{g' (X)}. \label{GnRHlayer2X}
\end{eqnarray}
\end{subequations}
Hence we have obtained a two-dimensional slow-fast system with slow variable $X$ and fast variable $x$. \\

\noindent \textbf{Boundary-layer system during the transitions.}
During fast motions 2 and 4, $(X,Y)$ evolves according to the $X$ time scale and the slowest variable $Y$ is almost constant. By setting $\delta =0$ in the rescaled system
\begin{subequations} \label{BLS3D}
\begin{eqnarray}
\eps \dot{x} &=& -y+f(x), \label{BLS3Dx} \\
\dot{y} &=& a_0 x+a_1 y+a_2+c X, \label{BLS3Dy} \\
\dot{X} &=& -Y+g(X), \label{BLS3DX} \\
\dot{Y} &=& \delta \left( X+b_1Y+b_2 \right), \label{BLS3DY}
\end{eqnarray}
\end{subequations}
one obtains the so-called {\it Boundary-Layer System}
\begin{subequations} \label{Fast3D}
\begin{eqnarray}
\eps \dot{x} &=& -y+f(x), \label{Fast3Dx}\\
\dot{y} &=& a_0 x+a_1 y+a_2+c X, \label{Fast3Dy}\\
\dot{X} &=& -Y+g(X), \label{Fast3DX}
\end{eqnarray}
\end{subequations}
where $Y\simeq g(-\gamma)$ for fast motion 2 and $Y\simeq g(\gamma)$ for fast motion 4.

\subsection{Different dynamical regimes}\label{sec-regs}

\noindent \textbf{Surge.}
The surge corresponds to phase 3 of \S \ref{sec-sketch}, when $X$ passes from $X_{\max}$ to $\gamma$. The dynamics is governed by system \eqref{GnRHlayer2}. Initially $x$ decreases to reach the vicinity of the nullcline $X=\tilde f(x)$, where
\[
\tilde f(x)=-\frac{a_0 x+a_1 f(x)+a_2}{c}.
\]
This is coupled with a significant increase of $y$. Subsequently $x$ grows and $y$ decreases, moving at the rate given by the slowest time scale $O(1/(\delta\eps))$. \\

\begin{figure}[htbp]
\centering
\includegraphics[scale=0.65]{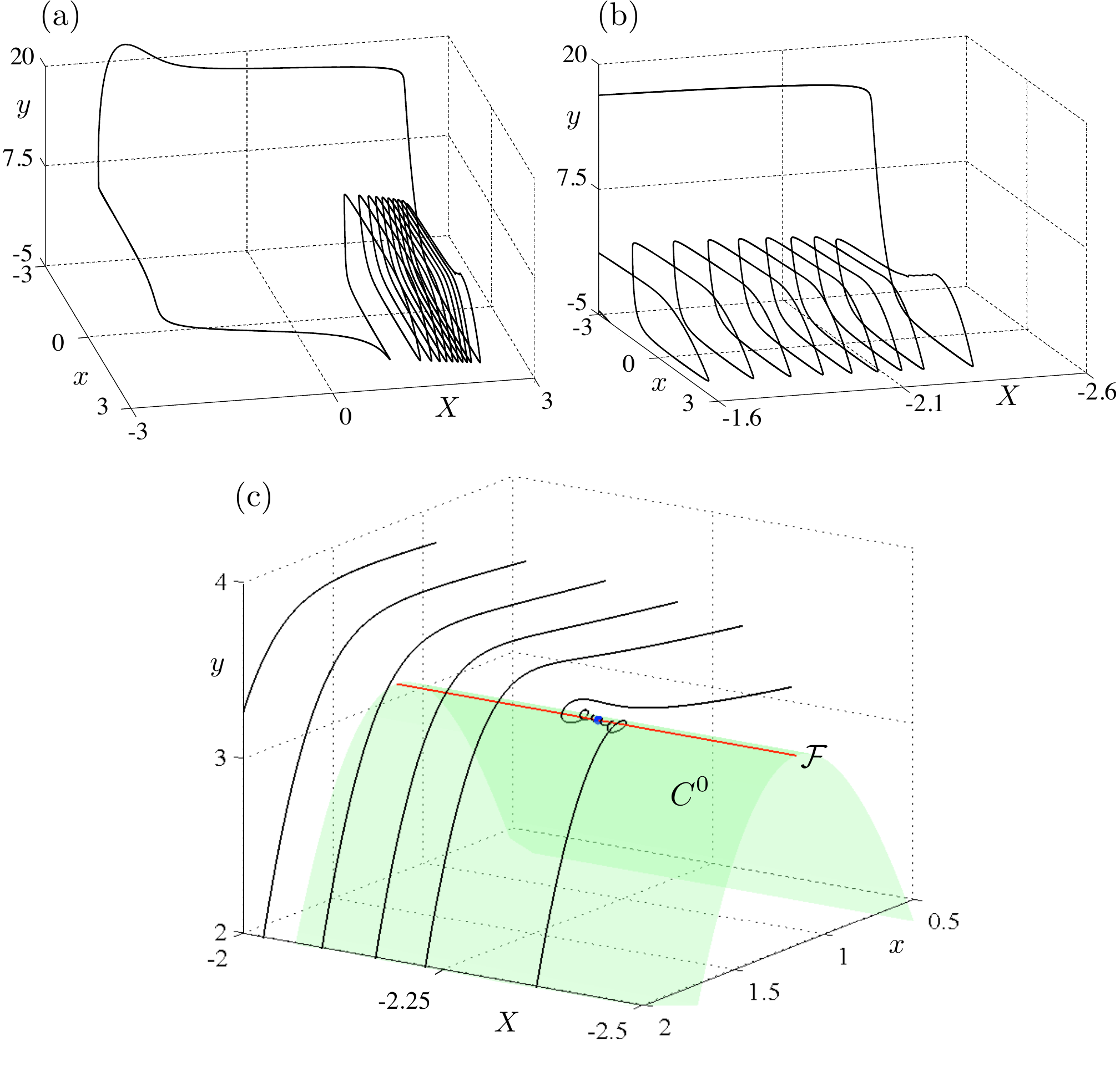}
\caption{Small oscillations corresponding to the post-surge pause. Panels (a) and (b) the full periodic orbit and a zoom in the region of the pause, respectively, in the three-dimensional phase space $(x,y,X)$; panel (c) shows a stronger zoom on the small oscillations of the pause, together with the critical manifold $C^0$, the fold curve $\mathcal{F}$ (red line) and the folded node (dot) of the three-dimensional subsystem \eref{3D3TS}.}
\label{mmo_3D_4D}
\end{figure}

\noindent \textbf{Small oscillations during the post-surge pause.}
Hypothesis (H1) guarantees that the surge is followed by a sequence of small oscillations taking place near the fold $(x_f, f(x_f))$. As $X$ has already reached the vicinity of $X_{\min}$ the dynamics is governed by \eqref{3D3TS} and can be described as follows. After the surge, the trajectory is attracted to a stable quasi steady state of node or focus type on the right branch of $y=f(x)$. As $X$ increases from $X_{\min}$, the quasi steady state changes stability as a pair of complex eigenvalues passes through the imaginary axis, that is, a slow passage through a Hopf bifurcation occurs \cite{an_87, an_88}. This phenomenon constitutes a delayed transition from  surge to pulsatility. The analysis of the small oscillations (see Figure \ref{mmo_3D_4D}) constitutes a large part of this article.

Recall that $X_f$ (introduced in hypothesis (H1)) is defined as the value of $X$ for which the $y$-nullcline intersects the cubic $x$-nullcline at its right fold. We will assume that $X_{\min} \lesssim X_f$ for the following reason. If $X_{\min}$ is significantly less than $X_f$ then the quasi steady state is still a node at the moment when the trajectory is attracted to it. During the passage near the fold the quasi steady-state turns into a stable focus and subsequently becomes unstable; first it is an unstable focus and then an unstable node. However, the trajectory is extremely close to the quasi steady state when it is a focus and only gets repelled from it after it has changed from unstable focus to unstable node. Hence no small oscillations can be seen. On the other hand, if $X_{\min} \lesssim X_f$ then the aforementioned quasi steady state can always be a focus, initially stable and subsequently unstable.

An important aspect of the dynamics are canards. A \textit{canard segment} is a segment of a trajectory which initially stays close to a stable branch of $y=f(x)$ for a time of $O(1/\eps)$, subsequently passes through the fold region and finally remains near the middle branch of $y=f(x)$ for a time of $O(1/\eps)$. A canard is a trajectory containing a canard segment. When a system possesses a folded singularity, there can be canard trajectories with small oscillations. We then define a $k$-th secondary canard as the canard trajectory making $k$ small oscillations near the fold. As part of the analysis we show that secondary canards separate the trajectories with different numbers of small oscillations. \\

\noindent \textbf{Pulsatility.}
Pulsatility is a region of transient relaxation oscillation corresponding to phase 1 of \S \ref{sec-sketch}. It is a direct continuation of the pause and the governing system is still \eqref{3D3TS}. As the slow nullcline of \eqref{eqdotx0}-\eqref{eqdoty0} cuts through the middle branch of $y=f(x)$ the dynamics is purely of relaxation type. \\

\noindent \textbf{Transitions.}
The first transition corresponds to phase 4 of \S \ref{sec-sketch} and follows the surge. The governing system is \eqref{Fast3D} with $Y=g(\gamma)$. The variables $(x, y)$ first evolve on the intermediate time scale $O(1/\eps)$, following the nullcline $y=f(x)$, subsequently jump to the right branch of the nullcline $y=f(x)$ and then follow the right branch of  $y=f(x)$, evolving on the time scale $O(1/\eps)$ and arriving to the vicinity of $(x_f,y_f)$ as $X$ reaches the vicinity of $X_{\min}$.

The second transition corresponds to phase 2 of \S \ref{sec-sketch} and precedes the surge. The governing system is \eqref{Fast3D} with $Y=g(-\gamma)$. Before $X$ reaches the vicinity of $X_{\max}$, the point $(x,y)$ can move down along the left branch of $y=f(x)$ and then turn back up  the left branch of $y=f(x)$ or it can jump to the right branch of $y=f(x)$ making another pulse before the surge. There is a canard phenomenon associated with this behavior which can incur some expansion. Estimating this expansion is a part of our analysis.

\subsection{Statement of the main theorem}
Before we can state our main theorem, we need to make one additional assumption:
\[
{\rm {\bf (H5)}} \qquad \frac{c}{\delta} \int_{x_{\rm sing}(X_{\max})}^{x_{\rm sing}(\gamma)} \frac{(\tilde f'(x))^2g'(\tilde f(x))dx}{f'(x)(\tilde f(x)+b_1g(\tilde f(x))+b_2)} 
>\frac{2}{\eps}\int_{-x_f}^{x_f} \frac{(f'(x))^2dx}{a_0x+a_1f(x)+a_2-cX_{f}},
\]
where
\[
\tilde f(x)=-\frac{a_0 x+a_1 f(x)+a_2}{c}.
\]
Hypothesis (H5) guarantees that the return map around the cycle is contracting. 
\begin{theorem}\label{th-main}
Provided that $\delta$ and $\eps$ are sufficiently small and (H1)-(H5) hold, there exists a unique stable limit cycle consisting of a number of small oscillations, a number of pulses and one surge. Some exceptional limit cycles, existing only in exponentially small parameter regions, contain canard segments. All the limit cycles are fixed points of a single passage around the cycle of surge, pause and pulsatility. Varying a regular parameter can lead to a change in the number of pulses or small oscillations by means of a passage through a canard explosion. There are two canard explosions, one associated with the upper fold and one with the lower fold. A passage through the canard explosion at the upper fold yields a transformation of a small oscillation to a pulse or vice versa. The passage through the canard explosion at the lower fold leads to an addition or a subtraction of a pulse.
\end{theorem}

\section{Folded singularities of system \eqref{3D3TS}} \label{sec-mr}

Folded singularities are usually studied in systems with two slow variables, however in this paper we need to consider system \eqref{3D3TS}, which has three time scales. The usual approach for classifying folded singularities is to consider the desingularized reduced system, see \cite{ps-mw_01} for instance. Here we will mimic this procedure for  \eqref{3D3TS}.

\subsection{Nature of the folded singularity} \label{sec-su}

Since the fastest variable in  \eqref{3D3TS} is $x$, we set the left hand side of  \eqref{3D3TSx} to $0$, obtaining the constraint $y=f(x)$. Hence the critical manifold corresponding to the fastest time scale is the cubic surface $C^0=\{y=f(x)\}$. It displays two folds respectively for
\[
x=\pm x_f = \pm \sqrt{\frac{\lambda_1}{-3\lambda_3}}, \quad y=f(x).
\]
By applying the procedure described at the beginning of \S \ref{sec-red-sys}, one obtains the desingularized reduced system
\begin{subequations} \label{drs}\
\begin{eqnarray}
\dot{x} &=& -(a_0 x+a_1 f(x)+a_2+c X) \label{drs_x},\\  
\dot{X} &=& -\delta\left ( \frac{X+b_1 (g(X)+O(\delta)) +b_2}{g'(X)+O(\delta)}\right )f'(x)=\Theta _{\delta }(x,X). \label{drs_X}
\end{eqnarray}
\end{subequations}
Note that the slow flow \eref{3D3TSy}-\eref{3D3TSX} on $C^0$ has the same orbits as the desingularized reduced system \eref{drs}, however one has to reverse the orientation where $f'(x)>0$, that is, on the repelling sheet of the critical manifold $C^0$.

The equilibria of system~(\ref{drs}) on the fold curve $\mathcal{F}$, that is, the folded singularities of the 3D system, are given by 
\begin{subequations} \label{eqdrs}
\begin{eqnarray}
f'(x_f) &=& 0, \label{eqdrs_x}\\
X_f     &=& -\frac{1}{c}(a_0x_f+a_1y_f+a_2), \label{eqdrs_X}
\end{eqnarray}
\end{subequations}
with $y_f$ given by $y_f=f(x_f)$. 
From the expression of $f(x)$, we get
\begin{subequations} \label{fs}
\begin{eqnarray}
x_f &=& \pm\sqrt{\frac{\minusone\!\lambda_1}{3\lambda_3}} \label{fs_x} \\
X_f     &=& \frac{\minusone\!1}{c} \left(\pm\sqrt{\frac{\minusone\!\lambda_1}{3\lambda_3}}(a_0+\frac{2}{3}\lambda_1a_1)+a_2\right), \label{fs_X}
\end{eqnarray}
\end{subequations}
with
\[
y_f=x_f(\lambda_3x_f^2+\lambda_1) = \pm\frac{2}{3}\lambda_1\sqrt{\frac{\minusone\!\lambda_1}{3\lambda_3}}.
\]
Note that 
\[
\frac{\partial \Theta _{\delta }}{\partial X} (X_f, x_f)=0
\]
due to the factor $f'(x)$. In addition we have $f''(x_f)=6\lambda_3x_f$. Hence, the jacobian matrix $J^{\eref{drs}}$ of system (\ref{drs}) at $(x_f,X_f)$ reads
\begin{eqnarray}
J^{\eref{drs}}(x_f, X_f)=\begin{pmatrix}
    -a_0                     & -c  \\
     -6\delta\lambda_3x_f\left ( \frac{X_f+b_1 (g(X_f)+O(\delta)) +b_2}{g'(X_f)+O(\delta)}\right )& 0  \\
\end{pmatrix}.
\end{eqnarray}
The eigenvalues of the matrix $J^{\eref{drs}}$ are given by
\begin{eqnarray}
\xi_{\pm} &=& \frac{1}{2}\left(-a_0\pm\sqrt{a_0^2+24c\lambda_3x_f \delta\left ( \frac{X_f+b_1 (g(X_f)+O(\delta)) +b_2}{g'(X_f)+O(\delta)}\right )}\right). \label{eig_fs}
\end{eqnarray}
It follows that if
\begin{equation}\label{eq-express}
X_{\rm eval}=24c\lambda_3x_f \left ( \frac{X_f+b_1 g(X_f) +b_2}{g'(X_f)}\right )<0
\end{equation}
then, for small enough $\delta$, there are two real eigenvalues of the same sign, i.e. the folded singularity is a folded node. Evaluating~(\ref{eq-express}) for $c=0.69$,  $a_0=1$, $a_1=0.02$, $a_2=0.8$, $b_1=0$, $b_2=-0.8$, $\lambda_3=-1$, $\lambda_1=1.5$, $\mu_3=-1$, $\mu_1=4$, one obtains $X_{\rm eval}= -3.3248$.

\subsection{Local form near the folded singularity}

We translate the origin to $(x_f, f(x_f), X_f)$, with $(x_f,X_f)$ given by \eqref{eqdrs}, and rescale the $x$ and $y$ variables. We first set
\begin{eqnarray*}
\overline{x} &=& \alpha \left( x-x_f \right), \\
\overline{y} &=& \alpha \left( y-f(x_f) \right),
\end{eqnarray*}
where $\alpha = \sqrt{-3\lambda_1 \lambda_3}$. In these new coordinates, system \eref{3D3TS} reads
\begin{eqnarray*}
\eps \delta \dot{\overline{x}} &=& -\overline{y}-{\overline{x}}^2 - \frac{1}{3\lambda_1} \overline{x}^3, \\
\delta \dot{\overline{y}} &=& a_0 \overline{x} + a_1 \overline{y} + \alpha \left( a_0 x_f + a_1 f(x_f) + a_2 + cX \right), \\
\dot{X} &=& \frac{X+b_1 (g(X)+O(\delta)) +b_2}{g'(X)+O(\delta)}.
\end{eqnarray*}
Now we translate the variable $X$ by setting
$\overline{X}=X-X_f$
to obtain the following system
\begin{subequations} \label{BLS_tr}
\begin{eqnarray}
\eps \delta \dot{\overline{x}} &=& -\overline{y} - {\overline{x}}^2 - \frac{1}{3\lambda_1} \overline{x}^3, \label{BLS_tr_x}\\
\delta \dot{\overline{y}} &=& a_0 \overline{x} + a_1 \overline{y} + \alpha c\overline{X}, \label{BLS_tr_y}\\
\dot{\overline{X}} &=& \varphi + \psi \overline{X} + O(\overline{X}^2, \delta). \label{BLS_tr_X}
\end{eqnarray}
\end{subequations}
where
\[
\varphi = \frac{X_f+b_1 g(X_f)+b_2}{g^{\prime}(x_f)}.
\]
Finally we introduce new variables $(x,y,z)$ defined by $(z,-x,y)=(\overline{x},\overline{y},\overline{X})$ and rescale the time to obtain the system:
\begin{subequations}\label{eq-weq}
\begin{eqnarray}
\dot{x} &=& \alpha c y -a_0 z + O(x), \label{eq-weq-x} \\
\dot{y} &=& \delta \left(\varphi +\psi y+ O(y^2, \delta) \right), \label{eq-weq-y} \\
\eps \dot{z} &=& x + z^2 + O(z^3). \label{eq-weq-z}
\end{eqnarray}
\end{subequations}
For $\delta=1$, system \eqref{eq-weq} is analogous to the normal forms considered in \cite{eb_90} and \cite{ps-mw_01}. Assuming that $\eps $ and $\delta $ are small constants (perturbation parameters), this system has one fast, one slow and one super-slow variable. Note that $X_{\rm eval}= 24c\lambda_3x_f\varphi$, $\lambda_3 <0$ and $c, x_f>0$. Hence, if \eqref{eq-express} holds then $\varphi>0$ and the folded singularity is a folded node. The case $\varphi<0$, or equivalently $X_{\rm eval}>0$, gives a folded saddle for $\delta$ sufficiently small. Finally $\varphi=0$ corresponds to a folded saddle-node. Other types of folded singularities do not occur for $\delta$ close to $0$.

\subsection{Local analysis near the folded node: statement of the result}\label{sec-local}

After dropping the \ $\bar \ $ \ signs and rescaling time, system (\ref{BLS_tr}) reads:
\begin{subequations} \label{BLS_tr1}
\begin{eqnarray}
\eps \dot{x} &=& -y - {x}^2 - \frac{1}{3\lambda_1} x^3, 
\label{BLS_tr1x}\\
\dot{y} &=& a_0 x + a_1 y + \alpha cX, \label{BLS_tr1y}\\
\dot{X} &=& \delta (\varphi +\psi X + O(X^2, \delta)). \label{BLS_tr1X}
\end{eqnarray}
\end{subequations}

Let the section $\Sigma^{\rm in}$ defined by $y=-\rho^2$, where $\rho >0 $ is small but fixed.  Let $S_{a,\eps}$ be the attracting Fenichel slow manifold, perturbed from the critical manifold $y=x^2+x^3/(3\lambda_1)$, near the section $\Sigma^{\rm in}$. Note that $\Sigma^{\rm in}$ is a transverse section of the flow of \eqref{BLS_tr1} intersecting $S_{a,\eps}$ close to the fold but $O(1)$ away from it. In this section we focus on describing the dynamics starting in the curve $S_{a,\eps}\cap\Sigma^{\rm in}$. Each trajectory starting in $S_{a,\eps}\cap\Sigma^{\rm in}$ enters the neighborhood of the fold, makes a number of small rotations, and then exits the fold region. The number of small oscillations can be different for different trajectories. {\em Canards} can now be defined as the trajectories that go into the repelling slow manifold $S_{r,\eps}$ (see \S \ref{sec-regs} for an alternative definition). A {\em $k^{th}$ secondary canard} is a canard that makes $k$ small oscillations in the fold region and subsequently runs into $S_{r,\eps}$.

Suppose the number of the small rotations for two trajectories $(x,y,X)$ and $(\tilde x, \tilde y, \tilde X)$ is different. Then there exists a secondary canard with initial condition somewhere on the segment between $(x,y,X)$ and $(\tilde x, \tilde y, \tilde Y)$. This way we can define sectors of the same rotation, or simply {\em sectors of rotation}, as the segments of $S_{a,\eps}\cap\Sigma^{\rm in}$ between the consecutive canards. We now state the main theorem of this section. This theorem leads to a precise definition and description of the sectors of rotation.

\begin{theorem}\label{th-seca}
There exists a number $R>0$ such that, for every $0<\nu<R$ there exists a family of $k^{th}$ secondary canards with 
\begin{equation*}
\frac{\nu}{\delta}<k<\frac{R}{\delta}.
\end{equation*}
The canards with consecutive rotation numbers are next to each other. The distance between the consecutive canards measured in the section $\Sigma^{\rm in}$ is bounded below by $C_1\delta\sqrt\eps$ and above by $C_2\delta\sqrt\eps$,
where $C_1$ and $C_2$ are positive constants. 
\end{theorem}

\begin{cor}\label{cor-sector}
The $k^{th}$ sector of rotation, defined as the region between the $k^{th}$ and the $(k+1)^{st}$ secondary canard consists of points whose trajectories make $k$ rotations in the fold region.
\end{cor}
Our proof of Theorem \ref{th-seca} is based on the application of the blow-up method and builds on the results of \cite{mk-mw_10}. This section is organized as follows. In \ref{sec-blupintro} we introduce the blow-up and its charts. In  \ref{sec-dK1} we analyze the dynamics in the entry chart K1, which covers a region near the critical manifold from $\Sigma^{\rm in}$ to $O(\sqrt\eps)$ away from the fold. In \ref{sec-dK2} we analyze the dynamics in the central chart K2, which describes the region very close to the singularity. We describe the delayed Hopf bifurcation occurring in this chart relying on the results of \cite{mk-mw_10}. In \ref{sec-proof} we prove Theorem \ref{th-seca}, building on the results of Sections \ref{sec-dK1} and \ref{sec-dK2}.

\subsection{Blow-up}\label{sec-blupintro}

We use the following blow-up function:
\begin{equation}
\begin{array}{rrcl}
\Phi: & \mathbb{R}_+ \times S^4    & \rightarrow & \mathbb{R}^5, \\
         & (\overline{r},\overline{x},\overline{y},\overline{X},\overline{\eps }) & \rightarrow &
                        (\bar{r}\bar{x}, \bar{r}^2 \bar{y}, \bar{r}\bar{X}, \bar{r}^2 \bar{\eps })=(x,y,X,\eps ).
\end{array}
\label{BU}
\end{equation}

In the entry chart $\overline{y}=-1$, the blow-up (\ref{BU}) transforms variables $(x,y,X ,\eps )$ into new variables $(x_1, r_1, X_1, \eps_1)$ as follows:
\begin{equation}\label{BU-K1}
x=r_1x_1,\quad y=-r_1^2,\quad X= r_1X_1,\quad  \eps=r_1^2\eps_1.
\end{equation}
After the transformation of system \eqref{BLS_tr1} and omitting a factor of $r_1$ (which corresponds to a time rescaling), we obtain the following equations:
\begin{subequations} \label{eqK1}
\begin{eqnarray}
\dot{x}_1 &=& -\frac{1}{2}x_1\eps_1F(x_1, r_1, X_1)-(-1 + {x_1}^2 +r_1\frac{1}{3\lambda_1} x_1^3), \label{eqK1_x}\\
\dot{r}_1 &=& \frac{1}{2}r_1\eps_1F(x_1, r_1, X_1), \label{eqK1_r}\\
\dot{X_1} &=& -\frac{1}{2}X_1\eps_1F(x_1, r_1, \eps_1)+\eps_1\delta (\varphi + \psi r_1 X_1 + O(r_1^2 X_1^2, \delta)),  \label{eqK1_X}\\
\dot{\eps_1} &=& -\eps_1^2F(x_1, r_1, X_1),\label{eqK1_e}
\end{eqnarray}
\end{subequations}
where $F(x_1, r_1, X_1)=-a_0x_1+a_1r_1-\alpha c X_1$.

In the transition chart $\bar\eps =1$, the blow-up (\ref{BU}) corresponds to the change of variables $(x, y, X, \varphi ,\eps )$ into $(x_2, y_2, X_2, \varphi_2, r_2)$ defined by
\begin{equation}
x=r_2 x_2,\quad y=r_2^2 y_2,\quad X= r_2 X_2,\quad  \eps=r_2^2,
\end{equation}
that is equivalent to
\begin{equation}\label{BU-K2}
x=\sqrt{\eps } x_2,  \quad y=\eps  y_2, \quad X=\sqrt{\eps } X_2.
\end{equation}
After the transformation of system \eqref{BLS_tr1}, canceling a factor of $\sqrt{\eps}$  and canceling the equation $\dot{\eps }=0$,
we obtain the following equations:
\begin{subequations} \label{eqK2}
\begin{eqnarray}
\dot{x_2} &=& -y_2 - x_2^2 - \frac{\sqrt{\eps }}{3\lambda_1} x_2^3, \label{eqK2_x}\\
\dot{y_2} &=& a_0 x_2 + \alpha cX_2 + a_1 \sqrt{\eps } y_2, \label{eqK2_}\\
\dot{X_2} &=& \delta \left(\varphi+ \sqrt{\eps} \psi X_2 + O(\eps X_2^2, \delta)\right). \label{eqK2_X^2}
\end{eqnarray}
\end{subequations}
To prove Theorem \ref{th-seca}, we need to find trajectories connecting from $S_{a,\eps}$ to $S_{r,\eps}$. As it will become clear from our forthcoming analysis, there is a natural extension of $S_{r,\eps}$ to $K1$. Hence we will be following the dynamics from $K1$ to $K2$ and then back to $K1$. Consequently, we need to be able to transform $K1$ to $K2$ and vice versa on the overlap of the charts.  The transformations between $K1$ and $K2$ are given by
\begin{equation}\label{K1toK2}
x_2=\frac{x_1}{\sqrt{\eps_1}},\qquad y_2=-\frac{1}{\eps_1},\qquad X_2=\frac{X_1}{\sqrt{\eps_1}}, \qquad\mbox{($K1$ to $K2$)},
\end{equation}
and 
\begin{equation}\label{K2toK1}
x_1=\frac{x_2}{\sqrt{-y_2}}, \qquad \eps_1=-\frac{1}{y_2}, \qquad X_1=\frac{X_2}{\sqrt{-y_2}}, \qquad\mbox{($K2$ to $K1$)}.
\end{equation}

\subsection{Extending Fenichel theory in chart $K1$}\label{sec-dK1}

The key observation concerning the dynamics of \eqref{eqK1} is that there exist center manifolds defined, approximately, by $x_1\approx \pm 1$. These manifolds are the extensions of the Fenichel slow manifolds $S_{a,\eps}$ and $S_{r,\eps}$. Near $r_1=0$, they intersect the hyperplane $\eps _1=0$ according to the equation
\[
-1 + {x_1}^2 +r_1\frac{1}{3\lambda_1} x_1^3 = 0.
\]
The center manifold $CM_a$ corresponding to $x_1\approx 1$ is attracting and given, near $r_1=0$, $\eps_1=0$, by the development
\[
x_1=1-\frac{1}{6\lambda_1} r_1+ O(\eps_1, r_1^2).
\]
Similarly, there exists an unstable center manifold $CM_r$, corresponding to $x_1\approx -1$.

The restriction of the flow to $CM_a$, after canceling a factor of $\eps _1$, which amounts to a time rescaling, is given by
\begin{subequations} \label{eqK1cm}
\begin{eqnarray}
\dot{r}_1 &=& \frac{1}{2}r_1\tilde F(r_1, X_1,\eps_1), \label{eqK1cm_r}\\
\dot{X_1} &=& -\frac{1}{2}X_1\tilde F(r_1, X_1, \eps_1) + \delta (\varphi + \psi r_1 X_1 + O(r_1^2 X_1^2, \delta )),  \label{eqK1cm_X}\\
\dot{\eps_1} &=& -\eps_1\tilde F(r_1, X_1,\eps_1),\label{eqK1cm_e}
\end{eqnarray}
\end{subequations}
with 
\[
\tilde F(r_1, X_1,\eps_1)=-a_0+(a_1+\frac{a_0}{6\lambda_1})r_1-\alpha c X_1+O(\eps_1, r_1^2)
\]
being the restriction of $F$ to $CM_a$. 

Hyperplanes $r_1=0$ and $\eps_1=0$ are invariant for \eqref{eqK1cm}. As $\varphi >0$, this system admits two singular points lying in $r_1=\eps_1 =0$ and defined by their $X_1$ component:
\begin{equation}\label{eq-X1pm}
X^{\pm}_1 = \frac{-a_0 \pm \sqrt{a_0^2-8\delta \alpha c \varphi}}{2\alpha c}.
\end{equation}
Both values of $X_1$ are negative and, since $\delta$ is small,
\[
-\frac{a_0}{\alpha c} \lesssim X^-_1 < X^+_1 \lesssim 0.
\]
Both singular points are hyperbolic.

\begin{figure}[htb]
\centering
\includegraphics[scale=0.45]{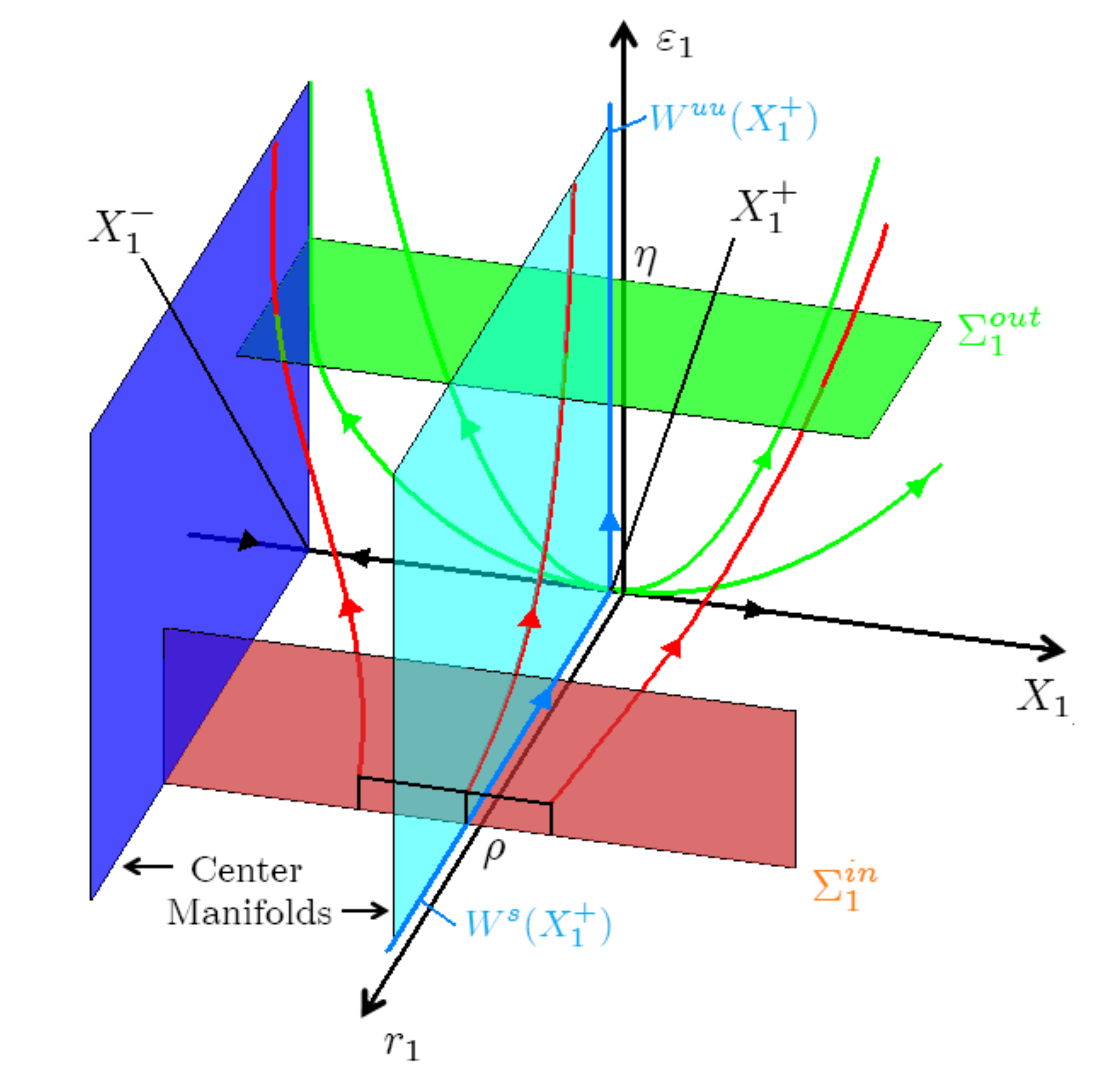}
\caption{Flow of \eqref{eqK1cm} obtained by reduction of the flow of the blown-up system in chart $K1$ to the attracting center manifold $CM_a$. It induces a transition between the sections $\Sigma^{\rm in}_1$ and $\Sigma^{\rm out}_1$ illustrate by the red orbits. The repulsive two-dimensional center manifold (sky blue) associated with the singular point $(0,X^+_1,0)$ (near the origin) contains the one-dimensional stable manifold $W^s(X^+_1)$ of this singular point the strong unstable manifold $W^{uu}(X^+_1)$. The attracting center manifold (purple) associated with the singular point $(0,X^-_1,0)$ is built in the same manner. Along the flow, the $X_1$ direction is expanded away from the sky blue manifold and contracted towards the purple one.}
\label{ChartK1}
\end{figure}

At each of the singular point $(r_1, X_1, \eps_1)=(0,X^{\pm}_1,0)$, the jacobian matrix $J_{\delta}^{\eref{eqK1cm}}$ associated with system \eref{eqK1cm} reads
\begin{multline*}
J_{\delta}^{\eref{eqK1cm}}(0,X^{\pm}_1,0) \\
=
\begin{pmatrix}
-\frac{1}{2}(a_0+\alpha c X^{\pm}_1) & 0 & 0 \\
-\frac{1}{2} X^{\pm}_1 \left( a_1+\frac{a_0}{6\lambda_1} \right) + \delta \left(\psi X^{\pm}_1 +O(\delta) \right) 
        & \frac{1}{2} \left(a_0+2\alpha c X^{\pm}_1 \right)  & -\frac{1}{2} O(1) \\
0 & 0 & a_0+\alpha c X^{\pm}_1
\end{pmatrix}
.
\end{multline*}
For $\delta \rightarrow 0^+$, one obtains at $(0,X^+_1,0)$,
\[
J_{0^{+}}^{\eref{eqK1cm}}(0,X^+_1,0)=
\begin{pmatrix}
-\frac{a_0}{2} & 0 & 0 \\
0^+  & \frac{a_0}{2}  & -\frac{1}{2} O(1) \\
0 & 0 & a_0
\end{pmatrix}
,
\]
and, at $(0,X^-_1,0)$,
\[
J_{0^{+}}^{\eref{eqK1cm}}(0,X^-_1,0)=
\begin{pmatrix}
0^- & 0 & 0 \\
-\frac{a_0}{2 \alpha c} \left( a_1+\frac{a_0}{6\lambda_1} \right)  & -\frac{a_0}{2}  & -\frac{1}{2} O(1) \\
0 & 0 & 0^+
\end{pmatrix}
.
\]
For $\delta>0$ small there exists a two-dimensional center manifold associated with the equilibrium $(0, X_1^-, 0)$. For $\delta=0$ this manifold is defined by the condition $\tilde F(r_1, X_1, \eps_1)=0$ and consists entirely of equilibria. For $\delta>0$ the flow on the center manifold becomes weakly hyperbolic. To estimate the corresponding eigenvalues, we use \eqref{eq-X1pm} to obtain the approximation
\begin{equation}\label{eq-X1pmapprox}
X_1^-=-\frac{a_0}{\alpha c}+\frac{2\varphi}{a_0}\delta+O(\delta^2).
\end{equation}
Using the fact that $a_0$, $\alpha$, $c$ and $\varphi$ are positive it is easy to see that the eigenvalue corresponding to the $r_1$-direction is negative and the eigenvalue corresponding to the $\eps_1$-direction is positive. Hence $r_1$ decreases and $\eps_1$ increases along trajectories. The flow of \eqref{eqK1cm}  is shown in Figure  \ref{ChartK1} and will be discussed in more detail below.

Based on the information collected above we can follow the passage of $CM_a$, and thus $S_{a,\eps}$, until the entry in the fold region. We define the following sections of the flow of \eqref{eqK1cm_r}-\eqref{eqK1cm_e}:
\[
\Sigma^{\rm in}_1=\{(x_1,r_1, X_1, \eps_1, \varphi_1)\;:\; r_1=\rho\},  
\]
where $\rho>0$ is the constant used in the definition of $\Sigma^{\rm in}$, and
\[
\Sigma^{\rm out}_1=\{(x_1,r_1, X_1, \eps_1, \varphi_1)\;:\; \eps_1=\eta\},  
\] 
where $\eta>0$ is a sufficiently small constant. We begin by explaining the meaning of $\Sigma^{\rm in}_1$ and  $\Sigma^{\rm out}_1$ in the context of the coordinates $(x, y, X)$ of system \eqref{BLS_tr1} and $(x_2, y_2, X_2)$ of system \eqref{eqK2}. The section  $\Sigma^{\rm in}_1$ corresponds to the section $\Sigma^{\rm in}$ defined at the beginning of \S \ref{sec-local}. Further it follows from \eqref{K1toK2} that the section $\Sigma^{\rm out}_1$ transforms to the section
\[
\Sigma^{\rm in}_2=\left\{(x_2, y_2, X_2) \;:\; y_2=-\frac{1}{\eta}\right\}. 
\]
Recall that $S_{a,\eps}$ is the attracting Fenichel slow manifold near the section $\Sigma^{\rm in}$. The set $S_{a,\eps}\cap\Sigma^{\rm in}$ is approximated by the line
\[
\left\{(x, y, X)\;:\; x=\rho, y=-\rho^2\right\}.
\]
To see which points in $\Sigma^{\rm out}_1$ can be reached by trajectories starting in $CM_a \cap \Sigma^{\rm in}_1$, we study the dynamics of system \eqref{eqK1cm} whose flow approximates the flow of \eref{eqK1} on $CM_a$. We see that the $X_1$-direction is expanded away from the equilibrium at the origin and is contracted towards the center manifold of the equilibrium $(0, X_1^-, 0)$. Hence the projection of the center manifold $CM_a$ onto the $X_1$-direction contains the interval $(-a_0/(\alpha c),\; 0)$. In fact the intersection of the center manifold with $\Sigma^{\rm out}_1$ is a thin band containing a line segment which is close to the interval
\[
x_1=1, \quad -\frac{a_0}{\alpha c} < X_1< 0.
\]

\subsection{Delayed Hopf bifurcation and the way in/way out function in chart $K2$}\label{sec-dK2}

We use system (\ref{eqK2}) to compute the way-in/way-out function near the fold. Note that \eqref{eqK2} is a slow-fast system with two fast and one slow variables, with singular parameter $\delta$. The critical manifold of \eqref{eqK2} is given by 
\[
x_2=-\frac{\alpha c}{a_0}X_2, \quad y_2=-\left (\frac{\alpha c}{a_0}\right )^2 X_2^2.
\]
We will denote this manifold by $S_0$.
The linearization of the fast system about $S_0$ is given by the matrix
\begin{equation}\label{matS0}
\left(
\begin{array}{cc}
2A(\eps)X & -1\\
a_0 & 0
\end{array}
\right) 
\end{equation}
with $A(\eps)=\alpha c/a_0+O(\sqrt\eps)$. Note that $A(\eps)>0$, for sufficiently small $\eps$. Hence, by Fenichel theory, given a constant $\eta > 0$,  there exist slow manifolds $S_{-,\delta}$ and $S_{+,\delta}$, attracting and repelling respectively, close to the line segments
\[
-\frac{1}{\eta}<X<-\eta \text{ and } \eta<X<\frac{1}{\eta}
\]
respectively.

\begin{figure}[htb]
\centering
\includegraphics[scale=0.45]{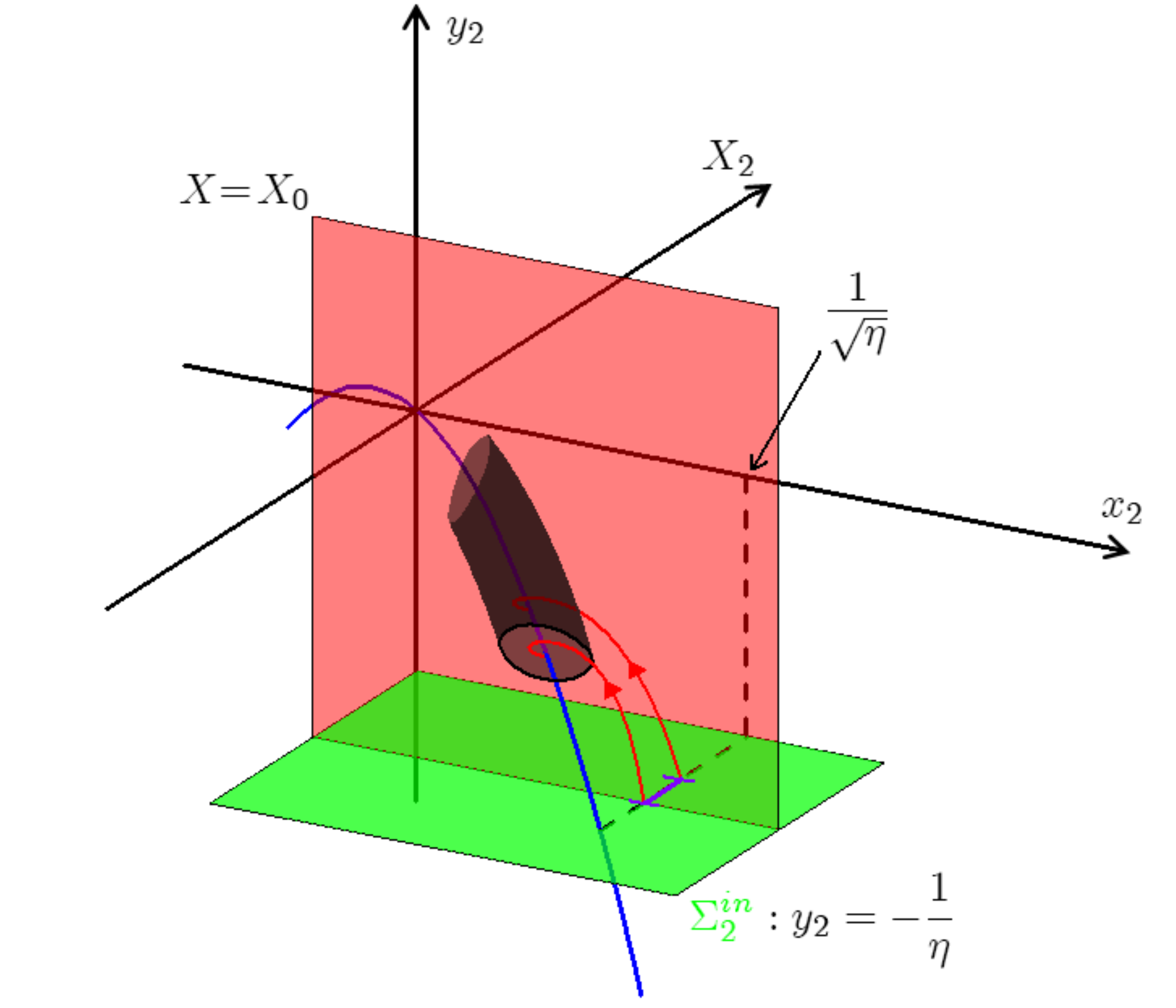}
\caption{Transition from $\Sigma_2^{\rm in}$ to the cylinder $C_{\delta}$ in chart $K2$. The image of $S_{a,\eps}\cap\Sigma^{\rm in}$ in $\Sigma_2^{\rm in}$ (purple interval) lies in $x_2=1/\sqrt{\eta}$. Under the flow of system \eref{eqK2}, trajectories starting from the purple interval reach the vicinity of the critical manifold $S_0$ by entering the cylinder $C_{\delta}$.}
\label{ChartK2Entry}
\end{figure}

\begin{figure}[htbp]
\centering
\includegraphics[scale=0.45]{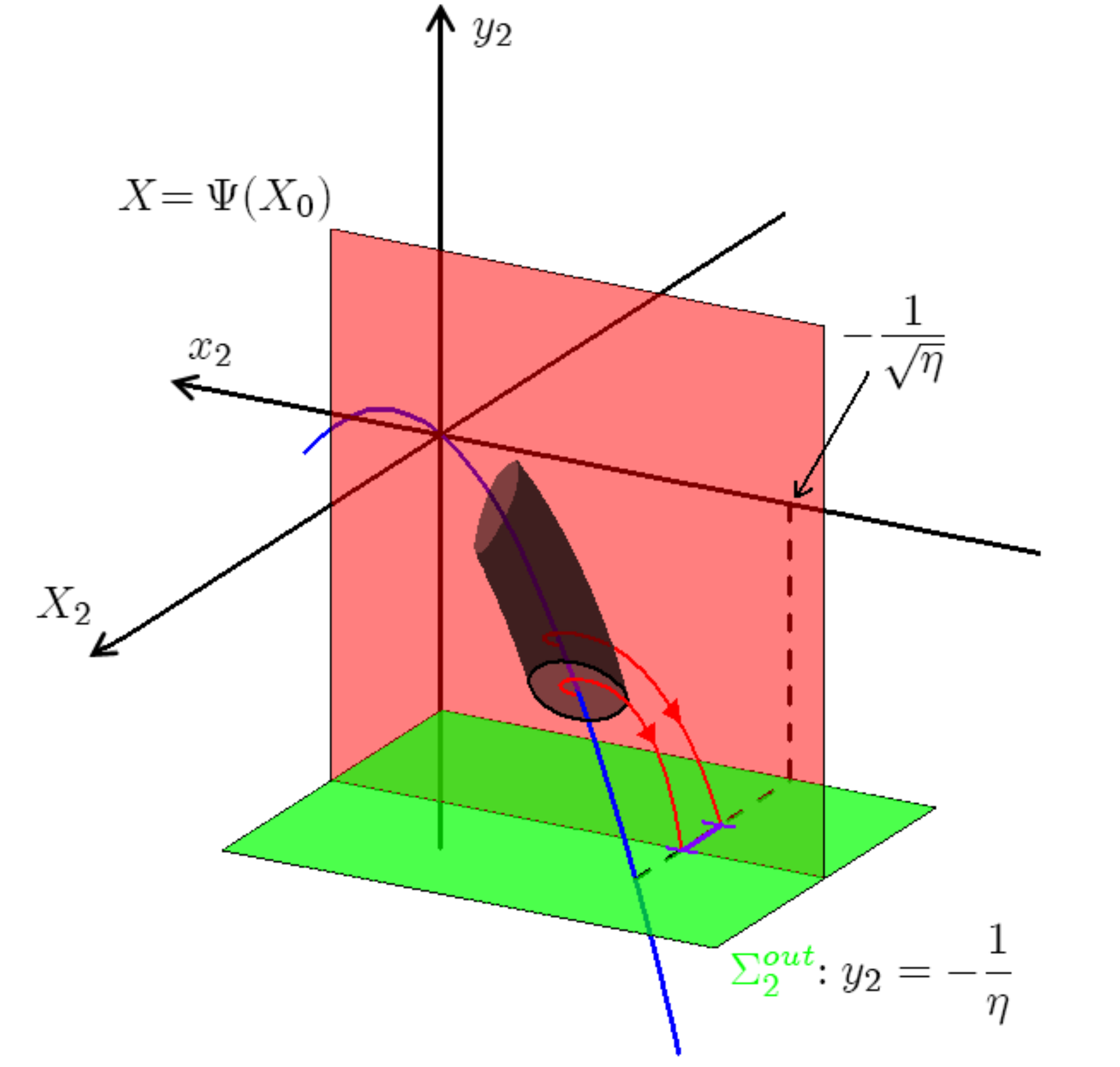}
\caption{ Transition from the cylinder $C_{\delta}$ in chart $K2$ to $S_{r,\eps}\cap\Sigma_2^{\rm in}$ (purple interval). The figure shows segments of trajectories of \eref{eqK2} 
starting in the interior of the cylinder $C_\delta$ and ending in $S_{r,\eps}\cap\Sigma_2^{\rm in}$. This figure can be understood as the `backwards in time' version of Figure \ref{ChartK2Entry}.}
\label{ChartK2Exit}
\end{figure}

We visualize the flow in the original coordinates $(x_2, y_2, X_2)$  with $X_2<0$, in Figure \ref{ChartK2Entry} ; the flow is strongly contracting in a tube around the critical manifold $S_0$. Similarly, Figure \ref{ChartK2Exit} shows the flow near $S_0$ with $X_2>0$; there the flow is strongly expanding.

We rectify $S_0$ by translating it to the line $(0,0,X_2)$, which is achieved by a transformation of the form: 
\[
x_2=\tilde x-\frac{\alpha c}{a_0} X_2+O(\sqrt\eps),\quad y_2=\tilde y-\left (\frac{\alpha c}{a_0}\right )^2 X_2^2+O(\sqrt\eps), \quad X_2=\tilde X.
\]
In the new variables, system \eref{eqK2} reads (after dropping the \ $\tilde .$ \ signs) 
\begin{subequations} \label{slowrect}
\begin{eqnarray}
\dot{x} &=& -y +2A(\eps)Xx -x^2+ O(\sqrt{\eps},\delta), \label{slowrect_x}\\
\dot{y} &=& a_0 x + O(\sqrt{\eps}, \delta, x^2), \label{slowrect_y}\\
\dot{X} &=& \delta(\varphi + \psi X + O(\sqrt{\eps } X^2)). \label{slowrect_X},
\end{eqnarray}
\end{subequations}
Note that $\sqrt\eps$ is a regular parameter in \eqref{slowrect}. Hence, to simplify the notation, we will suppress the dependance of $A$ on $\eps$. In Figure \ref{fig-cansec} we show the flow in the rectified coordinates for both $X$ negative and positive.

Note that the eigenvalues of the matrix given by \eqref{matS0} are off the real axis  if
\[
-\frac{\sqrt{a_0}}{A}< X < \frac{\sqrt{a_0}}{A}
\]
with negative real part for $X<0$ and with positive real part for $X>0$.
For $X<0$, we define the function $\Psi(X)$ by the formula
\begin{equation}\label{eq-psi}
\Psi(X)=X_*\quad\mbox{ with}\quad \int_{X}^{X^*}\frac{Z}{\varphi+\psi Z}\; dZ=0.
\end{equation}
The function $\Psi$ is the way in/way out function for all $X$ satisfying 
\[
-\frac{\sqrt{a_0}}{A}< X< 0.
\]
Heuristically this means that the trajectories attracted to $S_{-,\delta}$ near $X$, will be repelled from $S_{+,\delta}$ near $X_*$. To state a more precise result we introduce, for $X_0<0$, the sections
\[
\Sigma_{X_0}=\{(x,y,X)\;:\; X=X_0\}.
\]
The following result characterizes the transition map from  a section
$\Sigma_{X_0}$ to $\Sigma_{X_*}$, $X_*=\Psi(X_0)$.
\begin{prop}\label{prop-wiwo}
There exist constants $\eta_1>0$ and $\eta_2>0$ with the following property. For  any
\[
X_0 \in \left] -\frac{\sqrt{a_0}}{A},0\right[
\]
and $(x,y,X_0)\in\Sigma_{X_0}$ sufficiently close to the origin in $\Sigma_{X_0}$, let  $(x(t),y(t),X(t))$ be the trajectory of \eqref{slowrect} starting at $(x,y,X_0)$. Let $X_*=\Psi(X_0)$. There exists $t_*>0$ such that
\[
\begin{array}{c}
X(t_*)=X_*, \\
\forall t \in  ]0,t_*[,  \ \  X_0<X(t)<X_* \text{ and } x(t)^2+y(t)^2<\eta_1.
\end{array}
\]
Moreover, if the distance between $(x,y,X_0)$ and $S_{-,\delta}\cap \Sigma_{X_0}$  equals $\delta^{2(1+\alpha)}$ for some $0<\alpha<1/4$ then the distance between $(x(t_*),y(t_*), X_*)$ and $S_{+,\delta}$ is bounded below by $\eta_2\delta^{2(1+\alpha)}$.
\end{prop}

A proof of this result based on the work of Neishtadt \cite{an_88} can be found in \cite{mk-mw_10}, Corollary 5.1. Reference \cite{mk-mw_10} also contains extensions of Proposition \ref{prop-wiwo} to the case of $X_0< \sqrt{a_0}/A$, but we will not be concerned with these results in this article.

\subsection{Proof of Theorem \ref{th-seca}}\label{sec-proof}

In this section we will use both the original coordinates $(x_2, y_2, X_2)$ and the rectified coordinates $(x, y, X)$ in the following way. Let $C_\delta$ denote a cylinder of radius $\delta^{2(1+\alpha)}$ around $S_{-,\delta}$, where $0<\alpha<1/4$ is a constant. Trajectories starting in $\Sigma_2^{\rm in}$, with $X_2<0$ will enter $C_\delta$ (for small enough $\delta$) and subsequently exit $C_\delta$, with $X_2>0$. We will use the $(x_2, y_2, X_2)$ for the part of the trajectories outside $C_\delta$ and $(x, y, X)$ for the part of the trajectories inside $C_\delta$.

Recall that the interval
\[
x_1=1, \quad-\frac{a_0}{\alpha c} < X_1< 0
\]
is the image of $S_{a,\eps}\cap\Sigma^{\rm in}$ in $\Sigma^{\rm out}_1$. Translated to $\Sigma^{\rm in}_2$ this interval has the form
\[
x_2=\frac{1}{\sqrt{\eta}}, \quad -\frac{a_0}{\alpha c\sqrt{\eta}} < X_2< 0.
\]
Recall the center manifold $CM_a$ which coincided with the extension to $K_1$ of the slow manifold $S_{a,\eps}$. We denote the image of $CM_a$ in K2 by transformation \eqref{K1toK2} also by $CM_a$. Note that the point
\[
(x_2,y_2,X_2)=\left( \frac{1}{\sqrt{\eta}}, -\frac{1}{\eta}, -\frac{a_0}{\alpha c\sqrt{\eta}} \right)
\]
is on the critical manifold of system \eqref{eqK2}. By choosing $\eta<1/a_0$, we guarantee that the point where the eigenvalues change from real to complex, which  is given by
\[
X_2=-\frac{\sqrt{a_0}}{A},
\]
is included in the interval $(-a_0/(\alpha c\sqrt\eta), 0)$ provided that $\delta$ is small enough. Consider the segment of trajectory of \eqref{eqK2} starting at a point in $CM_a\cap \Sigma^{\rm in}_2$ such that
\[
-\frac{\sqrt a_0}{A} < X< -\eta_0<0
\]
and ending at a point in  $C_\delta$.  Note that the flow is predominantly in the fast directions and away from singularities. Hence the passage time is $O(\ln\delta)$ uniformly in $\delta$ and in $X_2\in (-\sqrt a_0/A,\, \eta_0)$. It follows that the $X_2$-coordinate of the endpoint of the segment of trajectory also satisfies
\[
-\frac{\sqrt a_0}{A} < X_2< \eta_0,
\]
provided that $\delta$ is small enough. 

We proceed using the rectified coordinates $(x,y,X)$. Consider a trajectory starting at a point in $C_\delta$ with $X=X_0$. Such trajectory first follows $S_{-,\delta}$ and subsequently $S_{+,\delta}$ until $X$ becomes approximately equal to $\Psi(X)$, see \eqref{eq-psi}. Corollary 5.1 in \cite{mk-mw_10} states that the flow  of \eqref{eqK2} for trajectories starting at $C_\delta$ is linear at lowest order (this occurs for the specified choice of $\alpha$). Moreover, for each $X_0$ verifying
\[
-\frac{\sqrt a_0}{A} < X_0< -\eta_0<0,
\]
there exists $\tilde X^*(X_0,\delta)$ such that the transition from  $C_\delta\cap \{X=X_0\}$ to $\{X=\tilde X^*(X_0,\delta)\}$ with
\[
\tilde X^*(X_0,\delta)=X^*(X_0)+O(\delta)
\]
is, at lowest order, a pure rotation by the angle $R(X_0)/\delta$, where
\[
R(X_0)=\int_{X_0}^{\Psi(X_0)}\frac{\sqrt{\frac{\alpha c Z}{a_0}-1}}{\varphi+\psi Z}\; dZ,
\]
We refer the reader to \cite{mk-mw_10}, \S 5, for further detail.

\begin{figure}[htbp]
\centering
\includegraphics[scale=0.5]{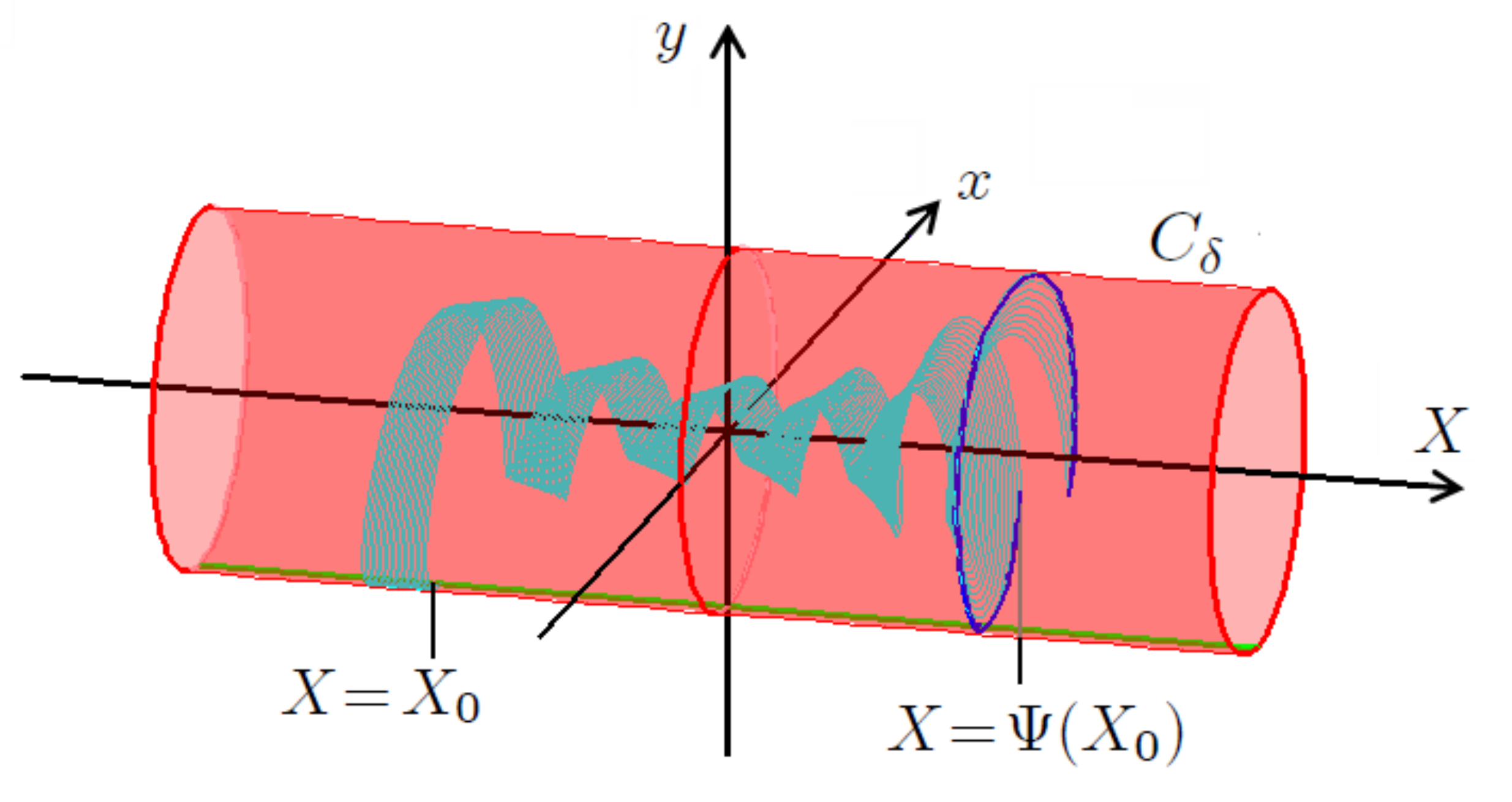}
\caption{Way-in/way-out transition near the rectified critical manifold of \eref{slowrect}. The trajectory entering the cylinder $C_{\delta }$ at $X=X_0<0$ is first attracted by the slow manifold $S_{-,\delta}$, remain near $S_{+,\delta}$ for a while and exits the cylinder at $X=\Psi(X_0)>0$. A $2\pi$ increment of the exit angle in $(x,y)$ is obtained by a $O(\delta)$ variation of the value $X_0$ of the entry.}
\label{fig-cansec}
\end{figure}

Consider an interval $(x_0,y_0, X)\in C_\delta$ with $X$ varying within $O(\delta)$ of $X_0$ and for every $X$ consider the image of $(x_0,y_0, X)$ by the transition via the flow of \eqref{eqK2} to the section $\{X=\tilde X^*(X_0,\delta)\}$. The image of the interval is a segment of a very tight spiral, as changing $X$ by $O(\delta)$ produces an increment of the angle greater than $2\pi$ while $\tilde X^*(X_0,\delta)$ changes very little. Now consider the continuation of $CM_r$ backwards in time, near $\tilde X^*(X_0,\delta)$. Since, backwards in time, the trajectories on $CM_r$ follow closely the fast fibers of \eqref{eqK2} and converge to $S_{+,\delta}$, they must intersect the mentioned spiral transversely. If, instead of taking the interval $\{(x_0,y_0, X)\in C_\delta , X=X_0+O(\delta)\}$, we take a segment of the continuation of $CM_a$ to $C_\delta$, we obtain a similar spiral and similar transverse intersections, separated by a distance bounded below by $K\delta$, for some constant $K>0$ (see Figure \ref{fig-cansec}).

These intersections correspond to secondary canards. A computation shows that
\[
-\frac{\sqrt a_0}{A}<X_0<-\eta  \Longrightarrow \Psi'(X_0)<0.
\]
Hence the number of rotations of the secondary canards monotonically increases as $X_0$ decreases. This implies that secondary canards are unique and that the consecutive canards differ by one rotation. It also follows that the distance between the secondary canards is $O(\delta)$. After translating back to the original coordinates $(x,y,X)$ (blowing down), we obtain the estimate of Theorem \ref{th-seca}.

\subsection{Passage through the fold region -- composition of the dynamics in the charts}

Using Theorem \ref{th-seca} we can now describe the dynamics in the fold region, starting from the union of the sectors of rotation in the section $\Sigma^{\rm in}$ to the exit from the fold region. We will restrict our attention to the trajectories on $S_{a,\eps}$ as all the other trajectories shadow a trajectory on $S_{a,\eps}$. We first consider the trajectories that are not close to canards. After passing through $C_\delta$ and exiting through its boundary for $X>0$ these trajectories separate quickly from $S_{r,\eps}$ and exit along the fast fibers. This is a simple and fast transition which does not incur much contraction. This transition is difficult to study mathematically due to resonance, but we will not focus on it here, referring the reader to \cite{mk-ps_01a}.

If the trajectories are close to a canard  they will arrive in $\Sigma^{\rm in}$ and must continue to either $\Sigma^{f}$, which is defined by $x_2=0$, (in the original coordinates by $x=x_f$), and subsequently reach $\Sigma^{\rm in}$, or they pass to the left branch of the nullcline resembling a {\it canard with head} and subsequently reach $\Sigma^{\rm in}$.

\section{Proof of Theorem \ref{th-main}} \label{sec-global}

\subsection{Contraction during the surge}\label{sec-contrasur}

\begin{figure}[htbp]
\centering
\includegraphics[scale=0.81]{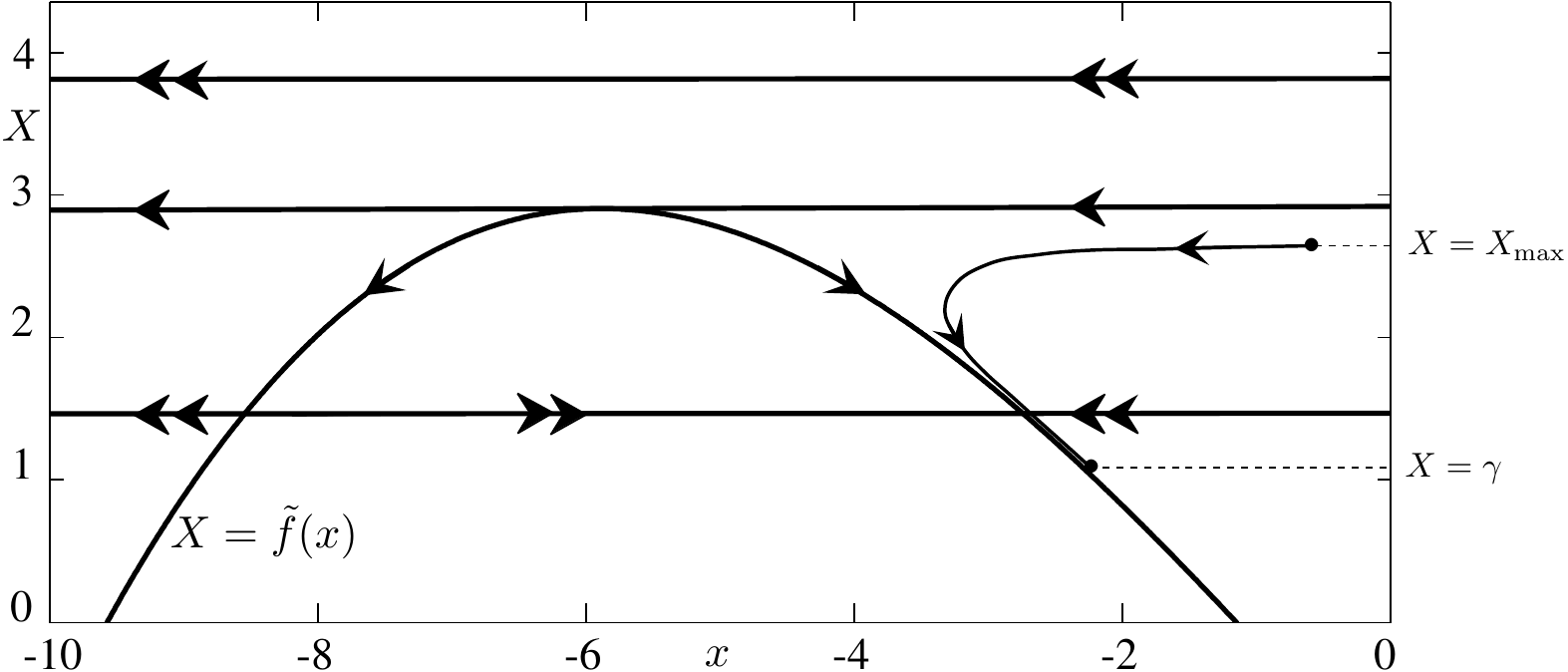}
\caption{Slow and fast dynamics of \eqref{GnRHlayer3} that approximate the flow of the whole system during the surge. Hypothesis (H3) ensure that $X$ is always smaller than the value of $X$ at the knee. Hence, at the beginning of the surge, $x$ is near $0$ and $X=X_{\max}$. Under the fast dynamics \eqref{GnRHlayer3x}, $(x,X)$ quickly reaches the slow manifold vicinity and then slowly goes down along it following the slow dynamics \eqref{GnRHlayer3X}. This mechanism is known to induce an exponential contraction between orbits.}
\label{pplanelayer2}
\end{figure}

Recall that surge begins as $X$ approaches $X_{\max}$ and the dynamics is governed by \eqref{GnRHlayer2}. We rewrite  \eqref{GnRHlayer2} for convenience using a different notation:
\begin{subequations} \label{GnRHlayer3}
\begin{eqnarray} 
\delta\dot{x} &=& c\frac{X-\tilde f(x)}{f' (x)}, \label{GnRHlayer3x}\\
\dot{X} &=& \frac{X+b_1g(X)+b_2}{g' (X)}, \label{GnRHlayer3X}
\end{eqnarray}
\end{subequations}
with
\[
\tilde f(x)=-\frac{a_0 x+a_1 f(x)+a_2}{c}=\frac{-a_1\lambda_3x^3-(a_0+a_1\lambda_1)x-a_2}{c}.
\]
Recall the definition of $x_{\rm sing}$ following the statement of Theorem \ref{th-main} and note that
\begin{eqnarray*}
\tilde f(x_{\rm sing}(X_{\max}))&=&X_{\max}, \\
\tilde f(x_{\rm sing}(\gamma))&=&\gamma.
\end{eqnarray*}
The slow manifold of \eref{GnRHlayer3} is defined by $X=\tilde f(x)$. It turns out that $\tilde f$ is an S shaped curve. To show that we compute the critical points. There are two, given by the formula:
\begin{equation}\label{crpts} 
x^2=-\frac{a_0+a_1\lambda_1}{3a_1\lambda_3}\approx 30.
\end{equation}
Further, $\tilde f''(x)=-6a_1\lambda_3x$, hence the negative critical point is a maximum and the positive one a minimum. Let $x_{c+}>0> x_{c-}$ be the critical points. Note that hypothesis (H3) is equivalent to the condition $X_{\max}< \tilde f(x_{c-})$. The slow and fast dynamics of \eqref{GnRHlayer3} are shown in Figure \ref{pplanelayer2}.

It follows from the slow-fast structure of \eqref{GnRHlayer3} and from the hypotheses (H3) and (H4) that the minimal (resp. maximal) value of $x$ during the surge is close to $x_{\rm sing}(X_{\max})$ (resp. $x_{\rm sing}(\gamma)$). The passage through surge is always an exponentially strong contraction, with contraction rate $O(\exp({\rm -C_3}/\delta))$, where $C_3>0$ is a constant. To compute $C_3$ we follow the approach of \cite{mk-ps_01b}, \S 2. We consider the reduced problem of \eqref{GnRHlayer3}, given by:
\begin{subequations} \label{GnRHlayerred3}
\begin{eqnarray} 
0 &=& -\tilde f(x)+ X, \label{GnRHlayerred3x}\\
\dot{X} &=& \frac{X+b_1g(X)+b_2}{g' (X)}, \label{GnRHlayerred3X}
\end{eqnarray}
\end{subequations}
or, parametrized by $x$,
\begin{equation}\label{GnRHlayerred3px}
x'=\frac{\tilde f(x)+b_1g(\tilde f(x))+b_2}{g' (\tilde f(x))\tilde f'(x)}.
\end{equation}
Let $x_0(t)$ be the solution of \eqref{GnRHlayerred3px}
defined on an interval $[0, t_{\rm endsurge}]$, with 
$x_0(0)=x_{\rm sing}(X_{\max})$ and $x_0(t_{\rm endsurge})=x_{\rm sing}(\gamma)$. Now, to estimate the contraction, we linearize the layer system
\begin{subequations}
\begin{eqnarray} 
\dot{x} &=& \frac{c(X-\tilde f(x))}{f' (x)}, \\
\dot{X} &=& 0.
\end{eqnarray}
\end{subequations}
Note that:
\[
\frac{\partial }{\partial x} \left[\frac{c(X-\tilde f(x))}{f' (x)} \right]_{|X=\tilde f(x)} = c\frac{\tilde f'(x)}{f'(x)}.
\]
Hence, the first order coefficient of the contraction rate is estimated by:
\[
C_3=c\int_{0}^{t_{\rm endsurge} }\frac{\tilde f'(x_0(t))}{f'(x_0(t))}dt.
\]
Changing the variables and using \eqref{GnRHlayerred3px} to express $dt$ as $dx/x'$ we get
\begin{equation}\label{eq-estC}
C_3=c\int_{x_{\rm sing}(X_{\max})}^{x_{\rm sing}(\gamma)} \frac{(\tilde f'(x))^2g'(\tilde f(x))dx}{f'(x)(\tilde f(x)+b_1g(\tilde f(x))+b_2)}.
\end{equation}

\subsection{Canard phenomenon during the passage from pulsatility to surge}\label{sec-passca}

While $X$ increases from $-\gamma$ to $X_{\max}$ the point $(x,y)$ may travel down along the left branch $y=f(x)$ towards the fold $(-x_f,-y_f)$. Subsequently it either turns back and travels towards $x_{\rm sing}(X_{\max})$ along the left branch $y=f(x)$ or jumps over to the other branch of $y=f(x)$ to complete another pulse. The transition between these two possibilities is a canard phenomenon, corresponding to the travel along the middle branch of $X$. We will not describe this canard phenomenon in detail, restricting our attention to the computation of the maximal expansion.

Proceeding as in \S \ref{sec-contrasur} we consider the slow flow of the system \eqref{Fast3D}, given by
\begin{subequations} \label{slowxX}
\begin{eqnarray} 
\dot{x} &=& \frac{a_0x+a_1f(x)+a_2+ cX}{f' (x)}, \label{slowxXx}\\
\dot{X} &=& -Y_0+g(X). \label{slowxXX}
\end{eqnarray}
\end{subequations}
and its solution $(x_0(t),X_0(t))$ with initial conditions $x(0)=-x_f$ and $X(0)=X_f$. The maximal amount of expansion occurs for trajectories traveling along the middle branch of $y=f(x)$ from the vicinity of $(-x_f,-y_f)$ to $(x_f, y_f)$ and is estimated by
\[
\int_{0}^{t_{\rm max} }f'(x_0(t))dt
\]
where $t_{\rm max}$ is the time when $x_0(t)$ reaches the upper fold $(x_f, y_f)$, i.e. $t_{\rm max}$ is defined by $x_0(t_{\rm max})=x_f$. Changing the variables, we get
\[
\int_{0}^{t_{\rm max} }f'(x_0(t))dt=\int_{-x_f}^{x_f} \frac{(f'(x))^2dx}{a_0x+a_1f(x)+a_2+ cX_0(t(x))}
\]
where $t(x)$ is defined by $x_0(t(x))=x$, with $-x_f\le x\le x_f$. Now since $-X_f<X(t)<X_{\min}$ we have
\[
\int_{-x_f}^{x_f} \frac{(f'(x))^2dx}{a_0x+a_1f(x)+a_2+ cX_0(t(x))} < 
\int_{-x_f}^{x_f} \frac{(f'(x))^2dx}{a_0x+a_1f(x)+a_2-cX_f}.
\]
Let
\begin{equation}\label{eq-estC4}
C_4=\int_{-x_f}^{x_f} \frac{(f'(x))^2dx}{a_0x+a_1f(x)+a_2-cX_{f}} .
\end{equation}
It follows that the amount of contraction incurred due to the canard phenomenon is approximately equal to $e^{C_4/\eps}$.

\subsection{Putting the pieces together}
In this section we put together all the phases of the dynamics, the surge, the pulsatility, the small oscillations, and the intermediate phases to get a transition around the entire cycle. We assume the $\delta=O(\eps)$, or more specifically, $\delta\le \eps$. We introduce four sections of the flow corresponding to the different phases of the dynamics.

\begin{figure}[htbp]
\begin{tabular}{cc}
\includegraphics[scale=0.41]{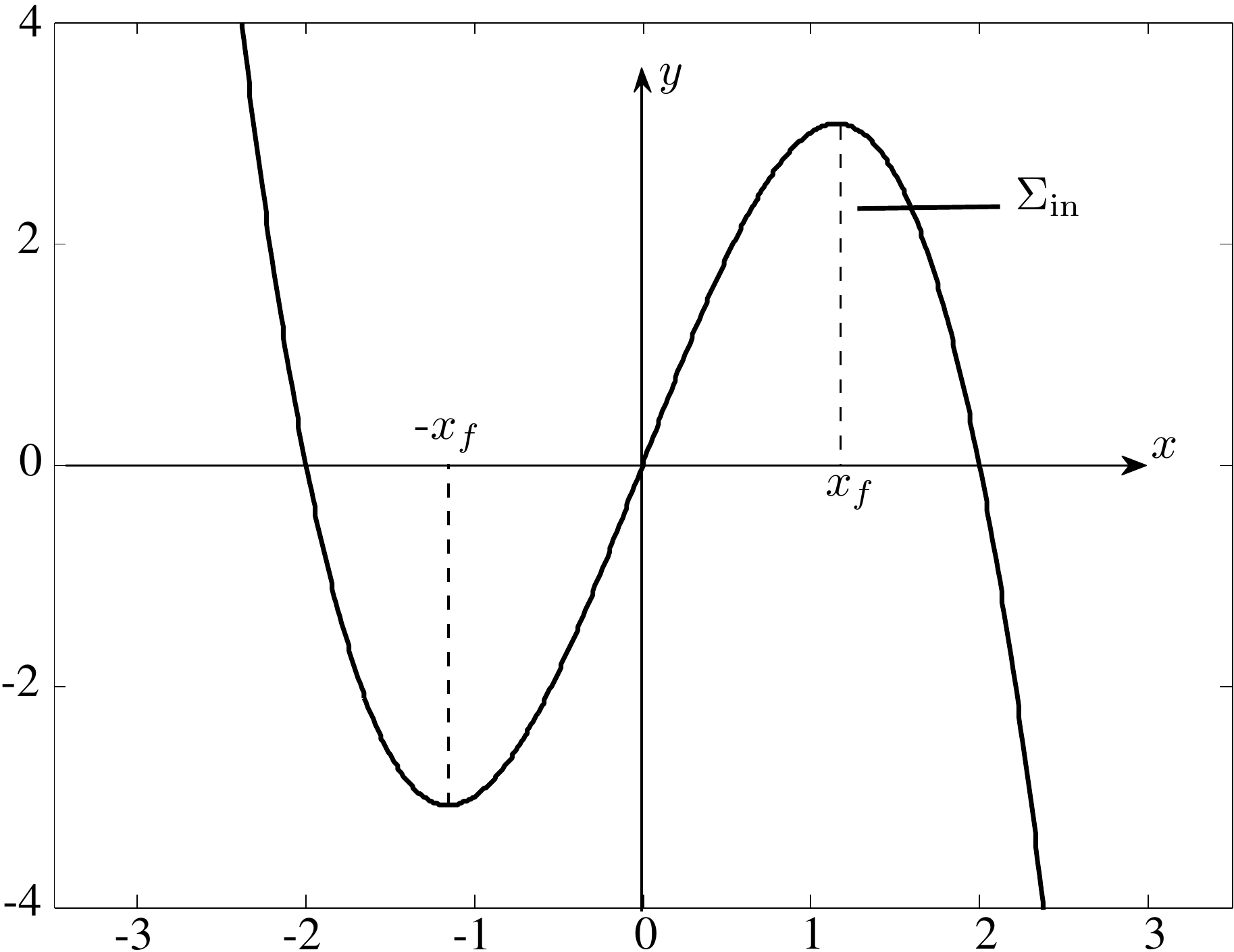} & \includegraphics[scale=0.41]{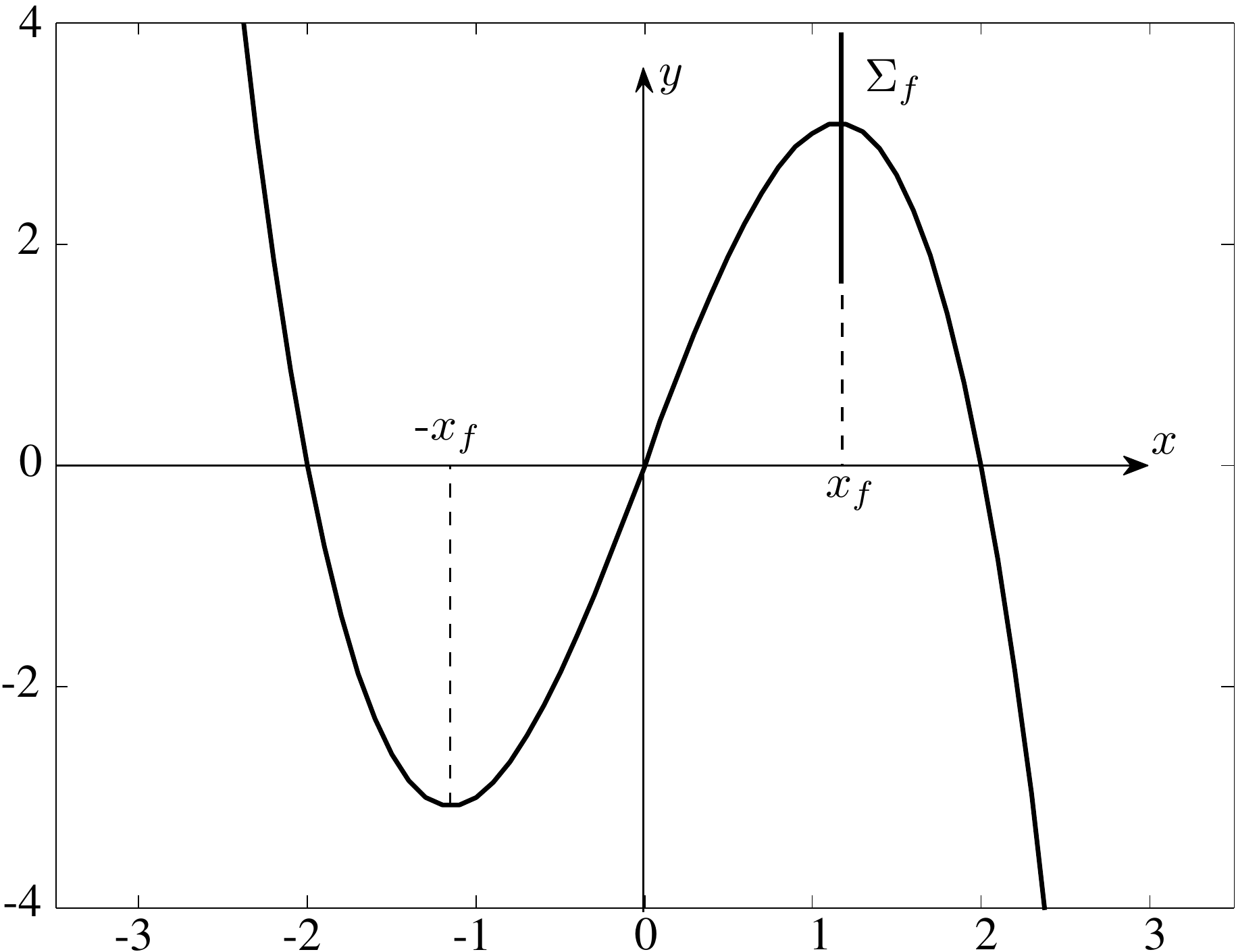} \\
\includegraphics[scale=0.41]{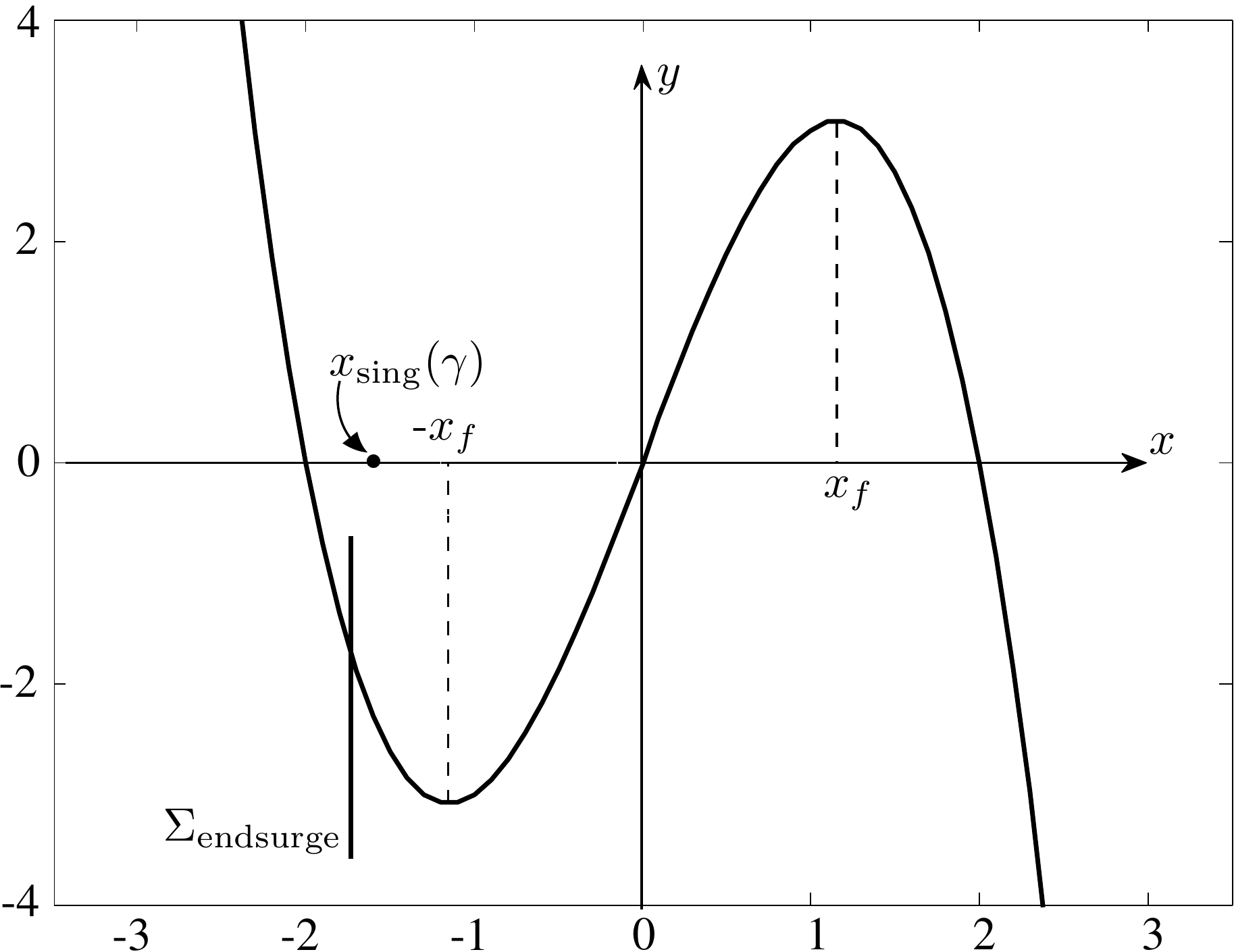} & \includegraphics[scale=0.41]{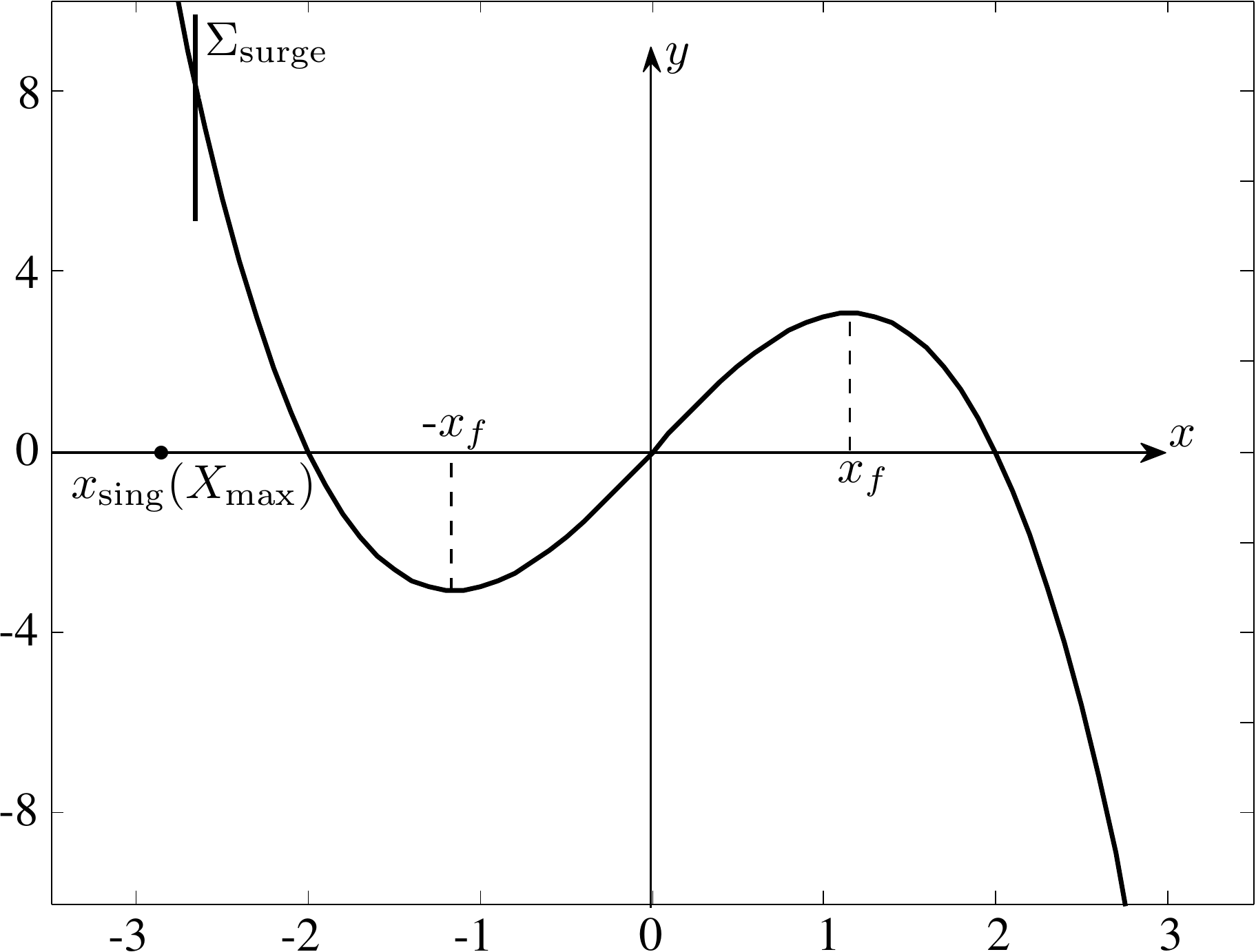}
\end{tabular}
\caption{The four sections of the flow of system \eqref{GnRHSystem} represented in the $(x,y)$-plane.}
\label{fig-sections}
\end{figure}

Let $\eta$ be a small constant. We can now define the sections of the flow:
\begin{align*}
\Sigma^{\rm in}&=\{(x,y,X,Y)\, :\, y=f(x_f)-\eta\}, \\
\Sigma^{f}&=\{(x,y,X,Y)\, :\, x=x_f\}, \\
\Sigma^{\rm surge}&=\{(x,y,X,Y)\, :\, x=x_{\rm sing}(X_{\max})+\eta\}, \\
\Sigma^{\rm endsurge}&=\{(x,y,X,Y)\, :\, x=x_{\rm sing}(\gamma)-\eta\}.
\end{align*}
The sections are shown in Figure \ref{fig-sections}. As shown in \S \ref{sec-contrasur} the transition map from $\Sigma^{\rm surge}$ to $\Sigma^{\rm endsurge}$ is a contraction with rate $O(e^{-C_3/\delta})$.

The transition from $\Sigma^{\rm endsurge}$ to $\Sigma^{\rm in}$ involves a passage through a folded node point and hence there are canard type phenomena occurring there. However, the passage of the Regulator through this point is fast ($O(1/\delta)$). Hence we can assume that the trajectories stay away from the canards so that no significant expansion is present.

If we restrict our attention to the interior of the sectors, staying away from the canards, then the transition from $\Sigma^{\rm in}$  to $\Sigma^{\rm surge}$ involves no significant expansion. Close passage to canards can be due to the folded node near the upper fold or due to the canard phenomenon during the passage from pulsatility to surge described in \S \ref{sec-passca}. 

We can now quickly describe what happens if the trajectories stay away from the canards. Consider the region in $\Sigma^{\rm in}$ consisting of the rotation sectors described in Theorem \ref{th-seca}. The image of this region in  $\Sigma^{\rm surge}$ is contained in a compact subset of  $\Sigma^{\rm surge}$. The image of this set by the transition from $\Sigma^{\rm surge}$ to $\Sigma^{\rm endsurge}$ is contained in a ball of radius $O(e^{-C_3/\delta})$, which we denote by $B_\delta$. By adjusting a parameter, for example $c$, we can arrange that the image of $B_\delta$ in $\Sigma^{\rm in}$ is contained in the union of the rotation sectors. We first consider the case when the image of $B_\delta$ in $\Sigma^{\rm in}$ is not close to a canard. We can now consider the return map from $\Sigma^{\rm endsurge}$ to itself restricted to $B_{\delta}$. This transformation is well defined (maps $B_\delta$ into itself) and is an exponential contraction. This proves the existence of a unique stable periodic orbit.

\begin{figure}[htbp]
\centering
\includegraphics[scale=0.65]{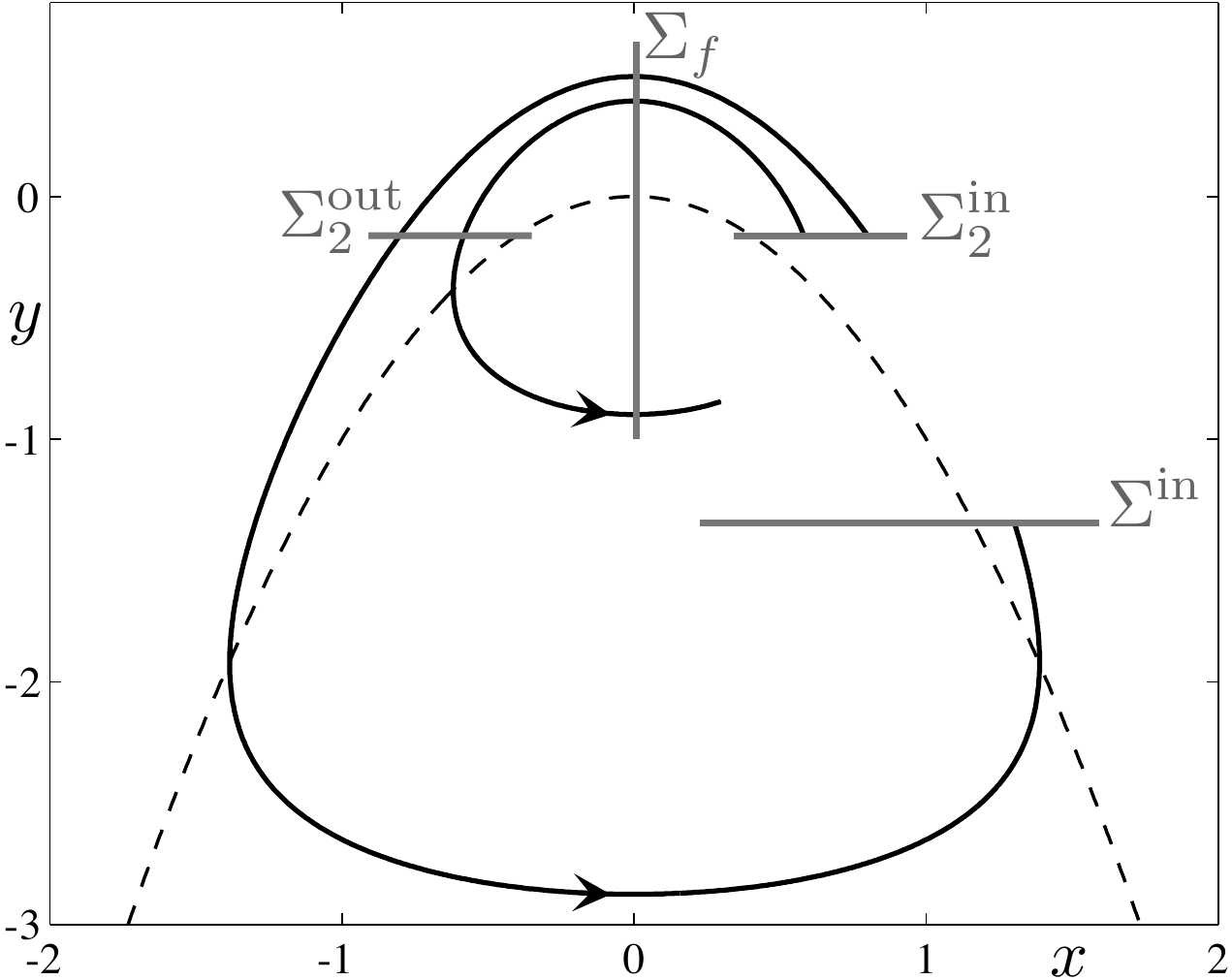}
\caption{A ``regular'' canard and a small canard.}\label{smaca}
\end{figure}

We now consider the case when trajectories starting in $\Sigma^{\rm surge}$ pass near canards. One way that a canard segment can be involved in the recurrent dynamics is if the image of $B_\delta$ in $\Sigma^{\rm in}$ intersects two sectors (with $\eps $ and $\delta$ small enough, it cannot intersect more due to the estimate on the size of the sectors given in Theorem \ref{th-seca}). There are now a few possibilities for how the trajectories can continue (see Figure \ref{smaca}).
\begin{enumerate}
\item Trajectories that are not close to a canard transition to the pulsatility stage in the same as in the case when no passage near canards was involved. Near such trajectories no significant expansion is incurred.
\item Trajectories that are close to canards, after passing through $K_2$, return to $\Sigma^{\rm in}$ through a segment of trajectory which looks like a canard cycle (with or without head). The maximal amount of expansion near such trajectories is incurred for maximal-like canards. The amount of expansion is estimated below.
\item Small canards, that return to $\Sigma^{f}$ and then to $\Sigma_2^{\rm in}$, without passing through $\Sigma^{\rm in}$. The amount of contraction near such trajectories is negligible.
\end{enumerate}
In both cases 2 and 3, $X-X_f$ is already positive as the trajectory reaches $\Sigma_2^{\rm out}$ so that a simple passage to pulsatility takes place following the return of the trajectory to $\Sigma^{in}$.

The maximal expansion by this transition is $e^{C_4/\eps}$, where $C_4$ is introduced in \eqref{eq-estC4}. To understand this estimate consider two points in $\Sigma^{\rm in}$ which are endpoints of two  trajectories starting in $\Sigma^{\rm in}$ at two points very close to each other but on the opposite side of the maximal canard. The flow backwards in time can contract the distance by the maximal contraction, backward in time, along the middle part of the fast nullcline, which is bounded below by  $e^{-C_4/\eps}$, with the constant $C_4$ computed analogously as $C_3$, see \S \ref{sec-contrasur}.

Another way that trajectories starting in $\Sigma^{\rm surge}$ may pass near a canard comes about by means of addition of a pulse at the end of pulsatility, or, in other words, by means of the canard phenomenon described in \S \ref{sec-passca}. As described in \S \ref{sec-passca} the maximal expansion is also bounded by $e^{C_4/\eps}$. The cumulative effect of the expansion coming from the two sources is $e^{2C_4/\eps}$.

Consequently, as long as
\begin{equation}
\frac{2C_4}{\eps} < \frac{C_3}{\delta},
\end{equation}
the return map from  $\Sigma^{\rm endsurge}$ to itself is an exponential contraction of $O(e^{-C_3/\delta+2C_4/\eps})$. Note that, for a fixed value of $\eps $, this condition is fulfilled for $\delta $ small enough. In the estimate, we include  the expansion incurred by the passage near both of the canard phenomena. Hence there exists a unique periodic orbit for every parameter value. This periodic orbit may contain a canard segment, which corresponds to a transition consisting of subtracting (respectively adding) a small oscillation and adding (respectively subtracting) a pulse at the beginning of the pulsatility stage, or it may include a canard segment corresponding to the addition of a pulse at the end of the pulsatility phase (the canard phenomenon described in \S \ref{sec-passca}).

\section{Numerical study}\label{sec-num}

We carried out numerical computation and continuation of periodic orbits in order to better describe the dynamics of our original four-dimensional system and visualize in more details the different transitions that shape the periodic orbits investigated in this work. We computed families of periodic orbits solutions of system \eqref{GnRHSystem} displaying pulses and surge, depending on various system parameters, using numerical continuation. We could then detect various canard-induced transitions affecting the number of pulses and the presence of a pause after the surge (see Figures \ref{pause_a2}, \ref{scenario_1} and \ref{scenario_2} below). We also computed attracting and repelling slow manifolds, as well as secondary canards, near the folded node of system~(\ref{BLS_tr1x})--(\ref{BLS_tr1X}), which approximates the behavior of the full system during the pause and explain its small oscillations (see Figure \ref{slowmancan} below). 

Systems with different time scales are well known to pose numerical problems because of their intrinsic stiffness, in particular when computing periodic orbits or, more generally, orbit segments \cite{jg-kh-ww_00,jg-ml_07}. The use of numerical continuation in the context of slow-fast dynamical systems has significantly increased over the last two decades. First, in the classic framework of limit cycle continuation, where very sharp transitions such as canard explosions could be finely rendered. Second, and more recently \cite{md-bk-ho_08}, in the context of manifold computation, where slow manifolds were approximated by families of orbit segments computed by continuing a parametrized family of two-point boundary value problems (BVP). The combination of orthogonal collocation to compute an orbit segment solution of a BVP, with a predictor-corrector algorithm to move one or several parametrized conditions of the BVP, proved very efficient compared to shooting methods. Indeed, solutions of singularly perturbed ODEs display very sensitive dependence on initial conditions and parameter variations. In this context, orthogonal collocation gives a better approximation of such an orbit segment by distributing the error along the orbit instead of accumulating it at one end point as with shooting. Furthermore, the continuation algorithm gives a good rendering of the piece of manifold of interest, where orbit segments are distributed according to arc length, hence, accounting for the changes of local curvature of the manifold. In this way, one can integrate slow-fast ODEs with suitable boundary conditions using the BVP solver embedded in numerical continuation packages such as \textsc{Auto} \cite{auto}. We now illustrate the use of numerical continuation tools to investigate system~(\ref{GnRHSystem}).

\subsection{Continuation of periodic orbits in parameter $a_2$}

\begin{figure}[H]
\centering
\includegraphics[scale=0.72]{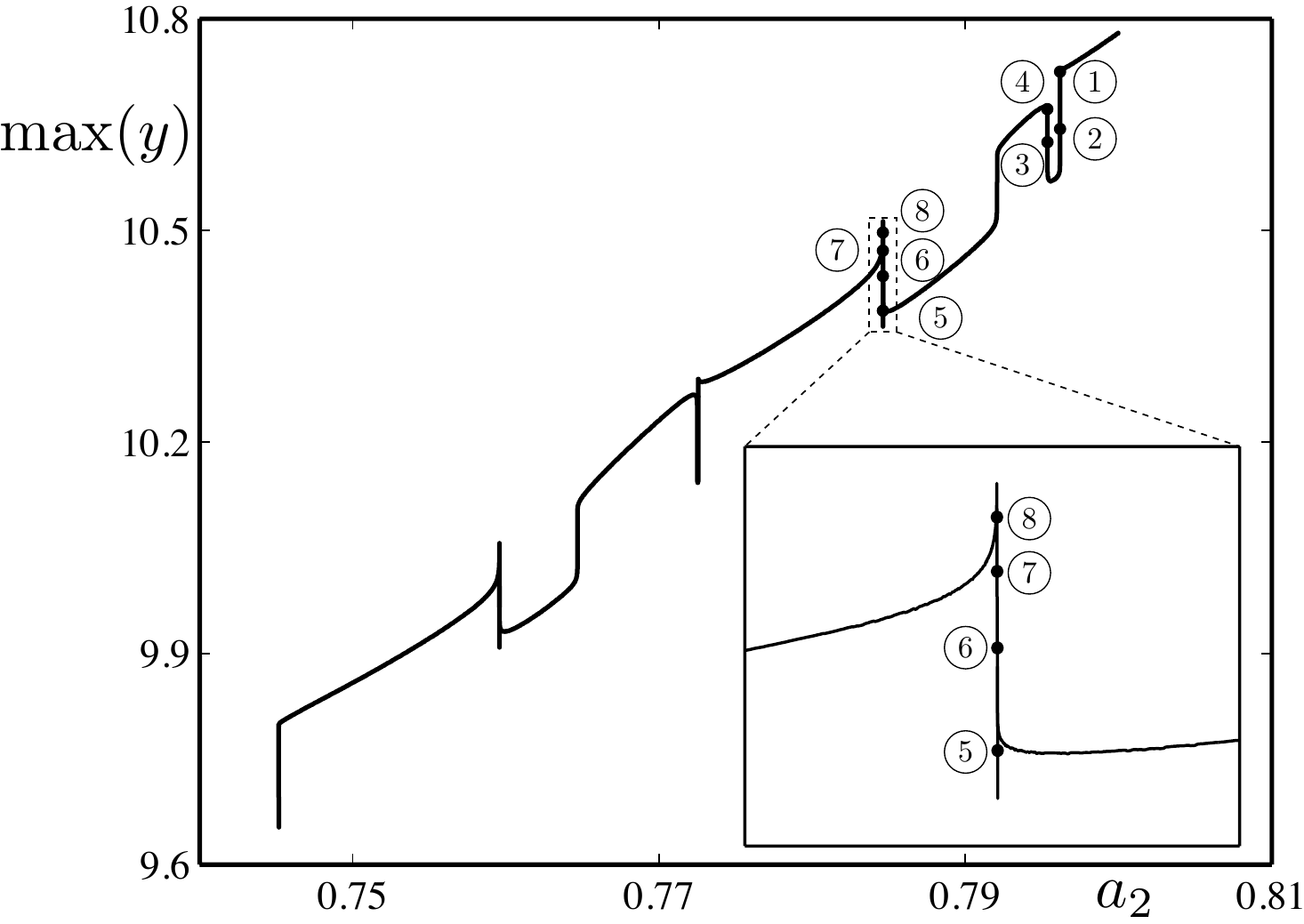}
\caption{Family of periodic orbits solution of the original system \eref{GnRHSystem} when $a_2$ is varied. The vertical axis shows the maximum in $y$ for each computed limit cycle along the branch. Eight orbits have been highlighted with black dots on the branch and given numbers from $1$ to $8$; they are shown in the two subsequent figures.}
\label{pause_a2}
\end{figure}

We start by a periodic orbit continuation that illustrates the various transitions, upon changes of parameter $a_2$, that take place in between different parts of the typical periodic orbit of system~(\ref{GnRHSystem}) as shown in Figure \ref{GnRHSecr}. The other parameters are fixed at values previously fixed, that is, 
$c=0.69$,  $a_0=1$, $a_1=0.02$, $a_2=0.8$, $b_1=0$, $b_2=-0.8$, $\lambda_3=-1$, $\lambda_1=1.5$, $\mu_3=-1$, $\mu_1=4$. The solution branch of periodic orbits is presented in Figure \ref{pause_a2} where we choose to display on the vertical axis the maximum in $y$ for each orbit as a measure of the solutions along the branch. The branch appears to be quite complicated with several rapid transitions that manifest themselves by quasi-vertical segments along the branch and that all have to do with canard trajectories. We identified two different types of transitions, affecting the periodic orbits at two different stages; before the surge, corresponding to the creation or annihilation of a pulse,  and after the surge during the pause, corresponding to the transition of a small oscillation to a pulse. These transitions correspond to the canard phenomenon discussed in \S \ref{sec-passca} and the canard phenomenon related to the small oscillations discussed in \S \ref{sec-su}. Each transition takes place within an exponential small variation of the parameter, therefore thus corresponding to a quasi-vertical segment on the branch. We will now describe our numerical results, which reveal an intricate sequence of canard explosions corresponding to the two types of transitions.

To be more specific, we introduce a labeling scheme for periodic orbits of \eqref{GnRHSystem}. We say that and orbit is of type $(p,s)$ if it involves $p$ pulses and $s$ small oscillations. Our first transition can be described as $(p,s) \to (p+1, s)$ (or  $(p,s) \to (p-1, s)$), while the second one as $(p,s) \to (p+1, s-1)$ ($(p,s) \to (p-1, s+1)$). We find two different scenarios in which both transitions occur. In the first scenario, they happen one after the other, which corresponds to two exponentially small bands of parameter values with associated quasi-vertical segment of the branch, separated by an order 1 interval of parameter values where the branch is, in comparison, quite flat. In the second scenario, they happen within the same exponentially small parameter variation. We isolate two sets of four orbits each along the branch, numbered $1$ to $4$ for the first set and $5$ to $8$ for the second, that undergo the first and the second transition, respectively. We now focus on each transition with the associated set of four chosen orbits. 

\begin{figure}[htbp]
\centering
\includegraphics[scale=0.42]{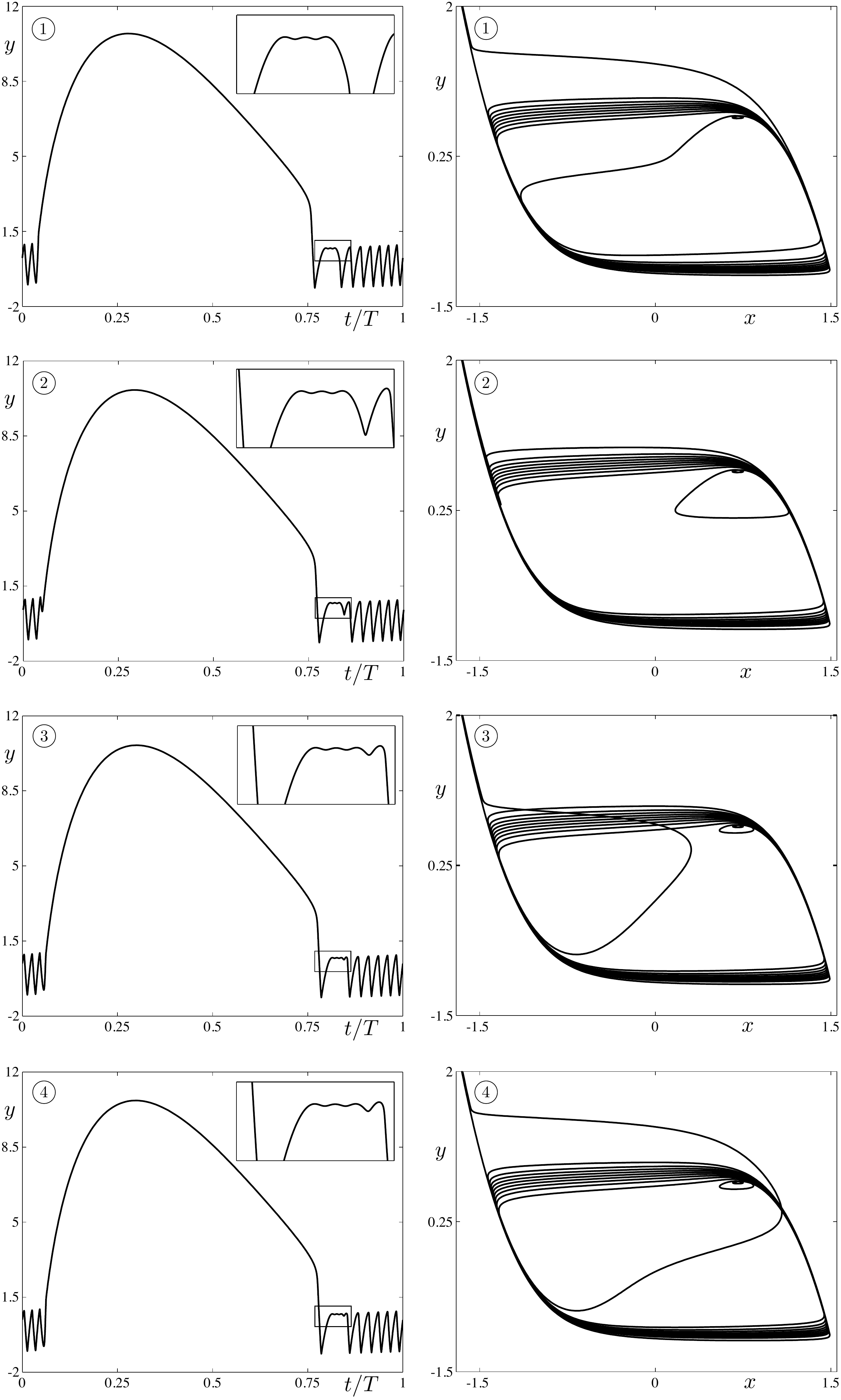}
\caption{First transition on the branch shown in Figure \ref{pause_a2} upon variation of $a_2$. From $1$ to $4$, the pause of the periodic attractor loses one small oscillations. We show the time profile of $y$ to illustrate this transition on the left panels and the projection of the orbit onto the $(x,y)$-plane on the right panels.}
\label{scenario_1}
\end{figure}

An example of the first scenario with both transitions is presented in Figure \ref{scenario_1}, where we show the time profile of $y$ for the orbits $1$ to $4$ on the branch in Figure~\ref{pause_a2}. This transition affects the number of small oscillations of the pause, this is why we enlarge each panel in the region of the pause and show the zoom in an inset; each panel is labeled with the number of the corresponding orbit in the solution branch. From orbit $1$ to orbit $4$, the pause gains one small oscillation, which corresponds to the loss of one pulse after the pause. However, one notices that there is very little difference between orbit $3$ and orbit $4$ in the pause. This is because this part corresponds to the second transition, where the change takes place before the surge with the appearance of one more pulse. In other words there are two canard explosions, which, using our labeling scheme, can be described as $(8,3) \to (7,4)$ and $(7,4)\to (8,4)$, giving a net result of an $(8,3) \to (8, 4)$ transition. Both transitions are canard-mediated, which one would see when plotting the orbits in the $(x,y)$-plane (see right panels of Figure~\ref{scenario_1}). However, within this first scenario they are separated by an $O(1)$ parameter interval.

\begin{figure}[!t]
\centering
\includegraphics[scale=0.41]{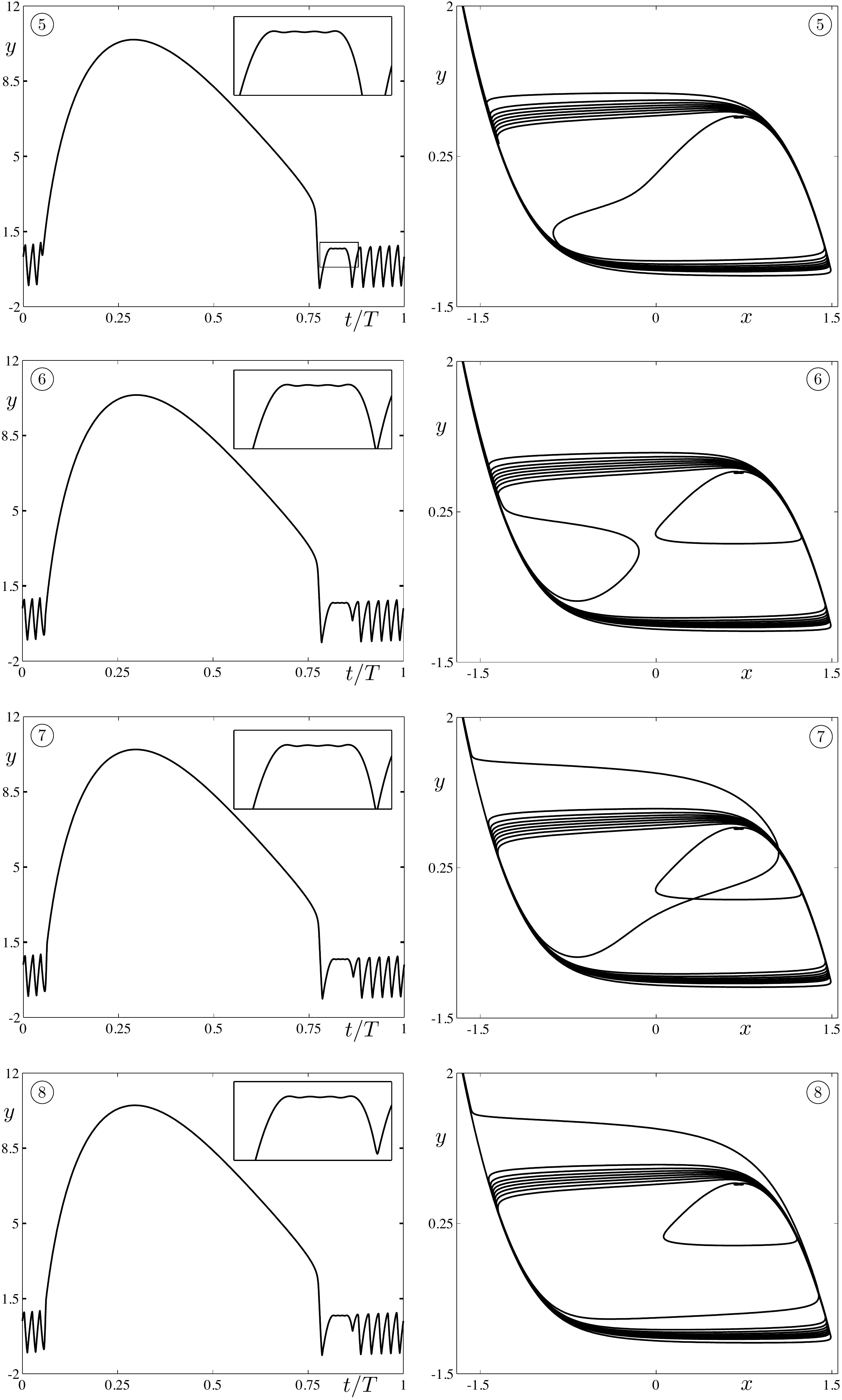}
\caption{Second transition on the branch shown in Figure \ref{pause_a2} upon variation of $a_2$. From $5$ to $8$, the start of the surge undergoes a canard explosion. We show the time profile of $y$ to illustrate this transition on the left panels and the projection of the orbit onto the $(x,y)$-plane on the right panels.}\label{scenario_2}
\end{figure}

\begin{figure}[!thb]
\centering
\includegraphics[scale=0.75]{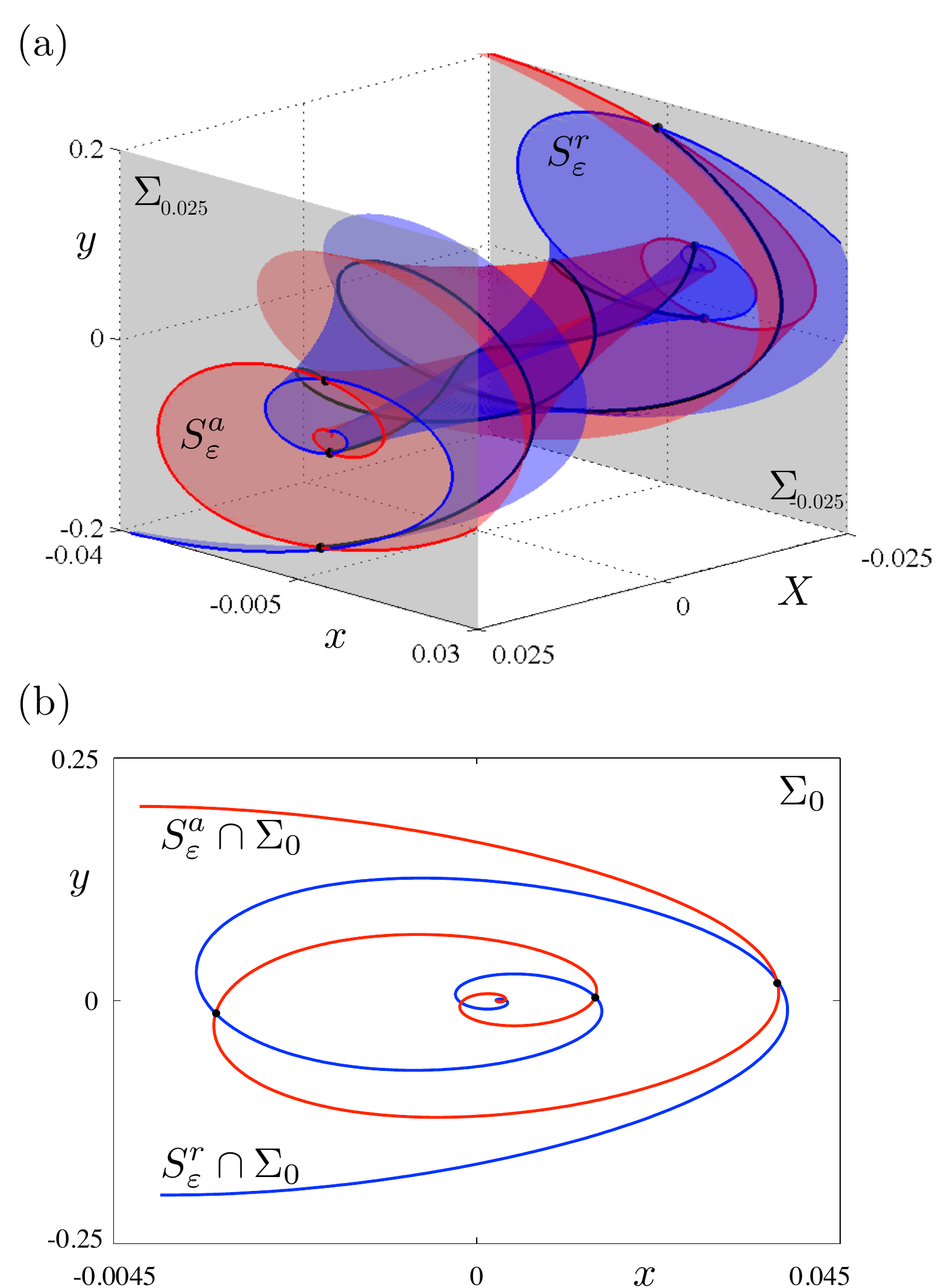}
\caption{Attracting ($S^a_{\eps}$) and repelling ($S^r_{\eps}$) slow manifolds of system (\ref{BLS_tr1}) near its folded node. Panel (a) shows a representation of these manifolds together with three secondary canards in the 3D phase space, in between the cross sections $\Sigma_{-0.025}:=\{X=-0.025\}$ and $\Sigma_{0.025}:=\{X=0.025\}$. Panel (b) shows the intersection curves of the slow manifolds $S^a_{\eps}$ and $S^r_{\eps}$ in $\Sigma_0:=\{X=0\}$.}\label{slowmancan}
\end{figure}

The second scenario is illustrated in Figure \ref{scenario_2} where orbits $5$ to $8$ from the solution branch are represented in the phase plane $(x,y)$. Here, both transitions seem to occur within the same exponentially small $a_2$-variation. Hence, one has two canard explosions in this plane : the first one corresponds to the transformation of a pulse into an additional small oscillation on the pause (visible essentially from orbit $5$ to $6$), the second one corresponds to the gain of one pulse before the surge (visible essentially from orbit $6$ to $8$). Hence, starting from $7$ pulses and $4$ small oscillations in case 5, one obtains $7$ pulses and $4$ small oscillations in case 8. The whole transition can be described as $(7,4) \to (6,5) \to (7,5)$. Note that the transitions in both scenarios are the same, the only difference is in the length of the parameter interval separating them. The two types of transitions are 
well explained by our theory, however we have no explanation for the very intriguing fact that,  in the second scenario, there are two canard explosions occurring simultaneously.

\subsection{Computation of slow manifolds and secondary canards on the pause}

We now illustrate the change of small oscillations on the pause, upon parameter variation, by computing slow manifolds and secondary canards of system~(\ref{BLS_tr1x})--(\ref{BLS_tr1X}), which represents a good approximation of the full system~\eref{GnRHSystem} in the region of the pause. This system is three-dimensional and possesses a folded node for the parameter values we consider (see section~\ref{sec-local}).
We computed slow manifolds and secondary canards near this folded node using the BVP strategy developed 
in~\cite{md-bk-ho_08,md-bk-ho_08_2}. That is, we approximate the manifolds by a one-parameter family of orbit segments with initial conditions moving on a curve traced on the critical manifold $C^0$, away from the fold $\mathcal{F}$, and end conditions restricted to a planar cross-section near the folded node. For simplicity, we take the one-dimensional manifold of initial conditions to be of the form $C^0\cap\{x=x_0\}$ (where $x_0$ is chosen so that this line is at a large enough distance from the fold curve), and the two-dimensional manifold in which the end conditions lie to be of the form $\Sigma_{\mathrm{end}}=\{X=X_{\mathrm{end}}\}$, with $X_0$ close to $0$ ($\Sigma_0$ corresponds to a cross-section containing the folded node). The dimensions of these manifolds of boundary conditions are chosen so that the resulting BVP is well-posed. 

In Figure~\ref{slowmancan}, we show the result of these manifold computations. Panel (a) shows an attracting slow manifold (red) $S^a_{\eps}$ and a repelling one (blue) $S^r_{\eps}$, computed in between sections $\Sigma_{-0.025}$ and $\Sigma_{0.025}$, together with three secondary canards (black curves) that correspond to transversal intersections between $S^a_{\eps}$ and $S^r_{\eps}$; we also show the intersection curves (red and blue curves) of the slow manifolds with both cross sections. The spiralling behavior of the slow manifolds is typical of the folded node scenario~\cite{mw_05,mb-mk-mw_06}. In panel (b), we present the intersection curves of the slow manifolds with the cross section $\Sigma_0$ that contains the folded node; once more, the figure is, as expected, very similar with previously computed slow manifolds in similar dynamical contexts~\cite{jg-kh-ww_00,mw_05,mb-mk-mw_06,md-bk-ho_08}.  

\section{Conclusion}\label{sec-dis}
In this paper we have studied the existence of MMOs in a system known as a phantom burster, consisting of two dimensional, unidirectionally coupled slow-fast oscillators or, in other words, one oscillator forcing the other by controlling the position of its null-cline. An additional feature of this system is the presence of three time scales; the dynamics of the forcing oscillator is slower than the dynamics of the one being forced. The orbits display alternatively three modes of dynamics: small oscillations near a fold, relaxation-type oscillations and a quasi steady state. Consequently, we needed to deal with quite complicated solutions, on the other hand the abundant structure of the system allowed us to obtain extensive results.

We have obtained two main results, one of global nature, namely the existence of a unique attracting periodic solution, and one of local nature, analyzing secondary canards of a folded node with an additional slow time scale. The local result relies strongly on the additional slow time scale and gives an elegant and rigorous proof of the existence of secondary canards and of the sectors of rotation. The global result relies on the local result and on the existence of strong contraction during the quasi steady-state phase of the dynamics. This result is rather elementary in nature, but it yields a surprising conclusion: in certain regions of the parameter space, the transition from an MMO with $n$ small oscillations to an MMO with $n+1$ small oscillations is free of complicated dynamics, a unique stable periodic orbit exists through the canard transition. This way we have obtained a complete characterization of the dynamics for all values of the control parameter. 

We were able to obtain rather strong results due to the simplicity of the system in question, but we hope that some of the ideas and techniques can be extended to other situations, involving three time scale dynamics and/or having the phantom burster structure. One simple but powerful idea we used was to identify the phase of the dynamics of the slowest system (fast, slow, quasi steady-state) to design a reduced system fitting the part of the dynamics in question. This way we derived important reductions that greatly simplified the analysis. Clearly this technique will be applicable in other contexts of multiple time scale systems, even if there is two way coupling between the systems with different time scales. 

The continuation results of Section \ref{sec-num} gave a very nice and accurate illustration of our theoretical results, as well as pointed to an intriguing phenomenon that we did not expect, namely a prediction of a double canard explosion. 
In this context it would be interesting to extend the continuation results to lower values of $\delta$, which is a numerical challenge. A more detailed numerical and theoretical study of the double canard explosion will be a subject of future work. 

We would like to point out that the system we have studied has been used to model different modes of GnRH (Gonadotropic Releasing Hormone) secretion and transitions between them. The biological mechanisms underlying the transitions between the surge mode and the pulsatility mode in the physiological GnRH secretion pattern are still poorly understood. Therefore, our study may contribute to the development of tools and insights that can be used in the study of this very important problem.

\section*{Acknowledgments}
This work has been financially supported by the large-scale initiative REGATE (REgulation of the GonAdoTropE axis) directed by Fr\'ed\'erique Cl\'ement: \\
\href{http://www.rocq.inria.fr/sisyphe/reglo/regate.html}{http://www.rocq.inria.fr/sisyphe/reglo/regate.html}. \\
The research of  M.K. has been funded by INRIA through a visiting professorship in the project-team SISYPHE and by the University of \'Evry-Val-d'Essonne through a visiting professorship in the Laboratoire Analyse et Probabilit\'es. \\
M.D. also acknowledges EPSRC through grant EP/E032249/1 and the Department of Engineering Mathematics at the University of Bristol (UK) where part of this work was completed. \\
The authors thank Jean-Pierre Fran\c{c}oise for helpful discussions.


\begin{thebibliography}{10}

\bibitem{eb_90}
E.~Beno{\^i}t.
\newblock Canards et enlacements.
\newblock {\em Publ. Math. IHES}, 72:63--91, 1990.

\bibitem{eb-jlc-fd-md_81}
E.~Beno{\^i}t, J.-L. Callot, F.~Diener, and M.~Diener.
\newblock Chasse au canard.
\newblock {\em Collect. Math.}, 31:37--119, 1981.

\bibitem{mb-mk-mw_06}
M.~{Br\o ns}, M.~Krupa, and M.~Wechselberger.
\newblock Mixed-mode oscillations due to generalized canard phenomenon.
\newblock In {\em Bifurcation Theory and spatio-temporal pattern formation},
  pages 39-- 64. Fields Institute Communications, 2006.

\bibitem{mb_88}
M.~Br{\o}ns.
\newblock Bifurcations and instabilities in the greitzer model for compressor
  system surge.
\newblock {\em Mathematical Engineering in Industry}, 2(1):51--63, 1988.

\bibitem{fc-jpf_07}
F.~Cl\'{e}ment and J.-P. Fran\c{c}oise.
\newblock Mathematical modeling of the {GnRH}-pulse and surge generator.
\newblock {\em SIAM J. Appl. Dyn. Syst.}, 6(2):441--456, 2007.

\bibitem{fc-av_09}
F.~Cl\'{e}ment and A.~Vidal.
\newblock Foliation-based parameter tuning in a model of the {GnRH} pulse and
  surge generator.
\newblock {\em SIAM J. Appl. Dyn. Syst.}, 8:1591--1631, 2009.

\bibitem{md-jg-bk-ck-ho-mw_12}
M.~Desroches, J.~Guckenheimer, B.~Krauskopf, C.~Kuehn, H.~M Osinga, and
  M.~Wechselberger.
\newblock Mixed-mode oscillations with multiple time scales.
\newblock {\em SIAM Review}, 54, 2012 (In press).

\bibitem{md-bk-ho_08}
M.~Desroches, B.~Krauskopf, and H.M. Osinga.
\newblock The geometry of slow manifolds near a folded node.
\newblock {\em SIAM J. Appl. Dyn. Syst.}, 7:1131--1162, 2008.

\bibitem{md-bk-ho_08_2}
M.~Desroches, B.~Krauskopf, and H.M. Osinga.
\newblock Mixed-mode oscillations and slow manifolds in the self-coupled
  fitzhugh-naguno system.
\newblock {\em Chaos}, 18:15107, 2008.

\bibitem{auto}
E.~J. Doedel, R.~C. Paffenroth, A.~R. Champneys, T.~F. Fairgrieve, Yu.~A.
  Kuznetsov, B.~E. Oldeman, B.~Sandstede, and X.~J. Wang.
\newblock Auto-07p: Continuation and bifurcation software for ordinary
  differential equations.
\newblock 2007.
\newblock Available at the URL: \texttt{http://indy.cs.concordia.ca/auto}.

\bibitem{be-mw_09}
B.~Ermentrout and M.~Wechselberger.
\newblock Canards, clusters and synchronization in a weakly coupled interneuron
  model.
\newblock {\em SIAM J. Appl. Dyn. Syst.}, 8(1):253--278, 2009.

\bibitem{jg_08}
J.~Guckenheimer.
\newblock Singular hopf bifurcation in systems with two slow variables.
\newblock {\em SIAM J. Appl. Dyn. Syst.}, 7:1355--1377, 2008.

\bibitem{jg-rh_05}
J.~Guckenheimer and R.~Haiduc.
\newblock Canards at folded nodes.
\newblock {\em Mosc. Math. J.}, 5:91--103, 2005.

\bibitem{jg-kh-ww_00}
J.~Guckenheimer, K.~Hoffman, and W.~Weckesser.
\newblock Numerical computation of canards.
\newblock {\em International Journal of Bifurcation and Chaos in Applied
  Sciences and Engineering}, 10(12):2669--2688, 2000.

\bibitem{jg-ml_07}
J.~Guckenheimer and M.~D. LaMar.
\newblock Periodic orbit continuation in multiple time scale systems.
\newblock In H.~M.~Osinga B.~Krauskopf and G.~Gal{\'a}n Vioque, editors, {\em
  Numerical Continuation Methods for Dynamical Systems}, pages 253--267.
  Springer, 2007.

\bibitem{mk-np-nk_08}
M.~Krupa, N.~Popovic, and N.~Kopell.
\newblock Mixed-mode oscillations in three time-scale systems: A prototypical
  example.
\newblock {\em SIAM J. Appl. Dyn. Syst.}, 7 (2):361--420, 2008.

\bibitem{mk-ps_01a}
M.~Krupa and P.~Szmolyan.
\newblock Extending geometric singular perturbation theoary to nonhyperbolic
  points -- fold and canard points in two dimensions.
\newblock {\em SIAM J. Math. Anal.}, 33:286--314, 2001.

\bibitem{mk-ps_01b}
M.~Krupa and P.~Szmolyan.
\newblock Relaxation oscillations and canard explosion.
\newblock {\em J. Differential Equations}, 174:312--368, 2001.

\bibitem{mk-mw_10}
M.~Krupa and M.~Wechselberger.
\newblock Local analysis near a folded saddle-node singularity.
\newblock {\em J. Differential Equations}, 248:2841--2888, 2010.

\bibitem{am-ps-hl-eg_98}
A.~Milik, P.~Szmolyan, H.~Loeffelmann, and E.~Groeller.
\newblock Geometry of mixed-mode oscillations in the 3d autocatalator.
\newblock {\em Int. J. Bifur. Chaos}, 8:505--519, 1998.

\bibitem{an_87}
A.~Neishtadt.
\newblock Prolongation of the loss of stability in the case of dynamic
  bifurcations {I}.
\newblock {\em Differ. Equ.}, 23:1385--1390, 1987.

\bibitem{an_88}
A.~Neishtadt.
\newblock Prolongation of the loss of stability in the case of dynamic
  bifurcations {II}.
\newblock {\em Differ. Equ.}, 24:171--176, 1988.

\bibitem{hr-mw-nk_08}
H.~Rotstein, M.~Wechselberger, and N.~Kopell.
\newblock Canard induced mixed-mode oscillations in a medial entorhinal cortex
  layer ii stellate cell model.
\newblock {\em SIAM J. Appl. Dyn. Syst.}, 7(4):1582--1611, 2008.

\bibitem{rw_07}
J.~Rubin and M.~Wechselberger.
\newblock Giant squid - hidden canard: the 3d geometry of the hodgkin huxley
  model.
\newblock {\em Biol. Cyber.}, 97(5), 2007.

\bibitem{ps-mw_01}
P.~Szmolyan and M.~Wechselberger.
\newblock Canards in $\mathbb{R}^3$.
\newblock {\em J. Differential Equations}, 177:419--453, 2001.

\bibitem{av-fc_10}
A.~Vidal and F.~Cl\'{e}ment.
\newblock A dynamical model for the control of the gnrh neurosecretory system.
\newblock {\em J. Neuroendocrinol.}, 22:1251--1266, 2010.

\bibitem{mw_05}
M.~Wechselberger.
\newblock Existence and bifurcations of canards in $\mathbb{R}^3$ in the case
  of the folded node.
\newblock {\em SIAM J. Appl. Dyn. Syst.}, 4:101--139, 2005.

\end{thebibliography}
\end{document}